\pdfoutput=1
\RequirePackage{silence}
\documentclass[a4paper,11pt]{amsart}
\usepackage[hmarginratio={1:1},vmarginratio={1:1},lmargin=60.0pt,tmargin=60.0pt]{geometry}

\allowdisplaybreaks

\usepackage[numbers]{natbib}
\usepackage[utf8]{inputenc}

\usepackage{latexsym,exscale,mathtools,textcomp}
\usepackage{amssymb,amsmath,amsthm,amsfonts,mathrsfs,bbm,enumitem,stmaryrd}
\usepackage[table]{xcolor}
\usepackage{graphicx}
\usepackage{mathtools}
\usepackage{ytableau}
\SetSymbolFont{stmry}{bold}{U}{stmry}{m}{n}

\usepackage{xparse}

\usepackage{dynkin-diagrams}

\setlist[enumerate]{itemsep=0.15cm,label=\emph{\upshape(\alph*)}}
\setlist[enumerate,2]{itemsep=0.15cm,label=\emph{\upshape(\roman*)}}

\usepackage{array}
\newcolumntype{C}{>{$}c<{$}}

\definecolor{mygray}{gray}{0.6}
\definecolor{mygraydark}{gray}{0.4}
\definecolor{mygraylight}{gray}{0.85}
\definecolor{spinach}{RGB}{46,139,87}
\definecolor{tomato}{RGB}{255,99,71}
\definecolor{orchid}{RGB}{143,40,194}
\definecolor{neon}{RGB}{77,77,255}
\definecolor{pumpkin}{RGB}{224,180,80}
\definecolor{citron}{RGB}{190,180,90}

\definecolor{lava}{RGB}{207,16,32}
\definecolor{cream}{RGB}{255,253,208}
\definecolor{verdigris}{RGB}{67,179,174}
\definecolor{Black}{RGB}{0,0,0}
\definecolor{mydarkblue}{RGB}{10,10,170}
\definecolor{darkspinach}{RGB}{20,70,20}
\definecolor{darktomato}{RGB}{155,40,30}
\definecolor{darkorchid}{RGB}{50,10,100}
\definecolor{darklava}{RGB}{150,8,16}

\usepackage{todonotes}

\let\emph\relax
\DeclareTextFontCommand{\emph}{\bfseries\em}

\newcommand{\placeholder}{{}_{-}}

\renewcommand{\dots}{\text{...}}

\let\<=\langle
\let\>=\rangle

\newcommand{\acts}{\centerdot}
\newcommand{\actsleft}{\mathop{\,\;\raisebox{1.7ex}{\rotatebox{-90}{$\circlearrowright$}}\;\,}}
\newcommand{\actsright}{\mathop{\,\;\raisebox{0ex}{\rotatebox{90}{$\circlearrowleft$}}\;\,}}
\renewcommand{\dots}{\text{...}}

\renewcommand{\ddots}{\raisebox{0.175cm}{\rotatebox{-45}{\text{...}}}}
\newcommand{\rddots}{\raisebox{-0.075cm}{\rotatebox{45}{\text{...}}}}

\DeclareMathOperator*{\medoplus}{\mathchoice
	{\textstyle\bigoplus}
	{\textstyle\bigoplus}
	{\scriptstyle\bigoplus}
	{\scriptscriptstyle\bigoplus}
}

\DeclareMathOperator*{\medotimes}{\mathchoice
	{\textstyle\bigotimes}
	{\textstyle\bigotimes}
	{\scriptstyle\bigotimes}
	{\scriptscriptstyle\bigotimes}
}

\DeclareMathOperator*{\medprod}{\mathchoice
	{\textstyle\prod}
	{\textstyle\prod}
	{\scriptstyle\prod}
	{\scriptscriptstyle\prod}
}

\DeclareMathOperator*{\medsum}{\mathchoice
	{\textstyle\sum}
	{\textstyle\sum}
	{\scriptstyle\sum}
	{\scriptscriptstyle\sum}
}

\newcommand{\medbinom}[2]{{\textstyle\binom{#1}{#2}}}

\DeclarePairedDelimiterX{\set}[1]{\{}{\}}{\setargs{#1}}
\NewDocumentCommand{\setargs}{>{\SplitArgument{1}{|}}m}{\setargsaux#1}
\NewDocumentCommand{\setargsaux}{mm}
{\IfNoValueTF{#2}{#1} {#1\,\delimsize|\,\mathopen{}#2}}%{#1\:;\:#2}

\newcommand{\ie}{\text{i.e.}}
\newcommand{\eg}{\text{e.g.}}
\newcommand{\cf}{\text{cf.}}
\newcommand{\etc}{\text{etc.}}
\newcommand{\aka}{\text{a.k.a.}}
\newcommand{\ver}{\text{verbatim}}

\newcommand{\muta}{\text{mutatis mutandis}}
\newcommand{\loccit}{\text{loc. cit.}}
\newcommand{\versus}{\text{vs.}}

\newcommand{\C}{\mathbb{C}}
\newcommand{\R}{\mathbb{R}}
\newcommand{\N}{\mathbb{Z}_{\geq 0}}
\newcommand{\Q}{\mathbb{Q}}
\newcommand{\Z}{\mathbb{Z}}
\newcommand{\K}{\mathbb{K}}

\newcommand{\ring}[1][R]{\mathbbm{#1}}

\newcommand{\bsym}[1][\lambda]{\boldsymbol{#1}}
\newcommand{\bsymv}[1][\mu]{\boldsymbol{#1}}
\newcommand{\stdvec}[1][i]{\bsym[\epsilon]_{#1}}

\newcommand{\vpar}{v}
\newcommand{\qpar}{q}
\newcommand{\vnum}[1]{[#1]_{\vpar}}
\newcommand{\vfac}[1]{[#1]_{\vpar}!}

\newcommand{\qnum}[1]{[#1]_{\qpar}}

\newcommand{\aform}{\mathbb{A}_{\vpar}}
\newcommand{\afrac}{\mathbb{A}_{\vpar}^{\prime}}
\newcommand{\apara}[1][{\bsymv}]{\mathbb{A}^{#1}_{\vpar}}
\newcommand{\afracpara}[1][{\bsymv}]{\mathbb{A}^{\prime,#1}_{\vpar}}

\newcommand{\para}{\lambda}
\newcommand{\parav}{\mu}

\newcommand{\Kpara}[1][{\bsym}]{\K_{\qpar}^{#1}}
\newcommand{\Kvpara}[1][{\bsym}]{\K_{\vpar}^{#1}}

\newcommand{\countingpara}{g}

\newcommand{\grdim}{\bsym[Z]\mathrm{dim}_{\K}}

\newcommand{\setstuff}[1]{\mathrm{#1}}

\newcommand{\obstuff}[1]{\mathtt{#1}}

\newcommand{\Aut}{\mathrm{Aut}}
\newcommand{\End}{\mathrm{End}}

\newcommand{\module}[1][M]{\obstuff{#1}}

\newcommand{\dgen}{\twoheadrightarrow_{d}}

\newcommand{\glgroup}[1][m]{\mathrm{GL}_{#1}(\C)}
\newcommand{\slgroup}[1][m]{\mathrm{SL}_{#1}(\C)}
\newcommand{\ogroup}[1][m]{\mathrm{O}_{#1}(\C)}
\newcommand{\sogroup}[1][3]{\mathrm{SO}_{#1}(\C)}
\newcommand{\spgroup}[1][m]{\mathrm{SP}_{#1}(\C)}

\newcommand{\sln}[1][n]{\mathfrak{sl}_{#1}}
\newcommand{\gl}[1][n]{\mathfrak{gl}_{#1}}
\newcommand{\ugl}[1][n]{U(\mathfrak{gl}_{#1})}
\newcommand{\uvgl}[1][n]{U_{\vpar}(\mathfrak{gl}_{#1})}

\newcommand{\uqgl}[1][n]{U_{\qpar}(\mathfrak{gl}_{#1})}
\newcommand{\uqsl}[1][n]{U_{\qpar}(\mathfrak{sl}_{#1})}
\newcommand{\uagl}[1][n]{U_{\mathbb{A}}(\mathfrak{gl}_{#1})}

\newcommand{\verma}[2][]{\obstuff{M}^{#2}}

\newcommand{\qverma}[2][\qpar]{\obstuff{M}^{#2}_{#1}}
\newcommand{\averma}[2][\mathbb{A}]{\obstuff{M}^{#2}_{#1}}

\newcommand{\amap}[2]{\phi^{#2}_{#1}}

\newcommand{\vsym}[1][{\bsym[k]}]{\obstuff{Sym}^{#1}_{\vpar}\K^{2}}
\newcommand{\qsym}[1][{\bsym[k]}]{\obstuff{Sym}^{#1}_{\qpar}\K^{2}}

\newcommand{\gt}[1][\lambda]{GT_{#1}}

\newcommand{\dense}[2][]{\obstuff{D}^{#2}}
\newcommand{\densemult}[2][]{\mathfrak{D}^{#2}}
\newcommand{\qdensemult}[2][]{\mathfrak{D}^{#2}_{\qpar}}

\newcommand{\qdense}[2][\qpar]{\obstuff{D}^{#2}_{#1}}

\newcommand{\sums}[1]{\scalebox{0.7}{$\Sigma$}{#1}}
\newcommand{\prods}[1]{\scalebox{0.7}{$\Pi$}{#1}}

\newcommand{\polyalg}[1][{\bsym}]{\obstuff{P}^{#1}}

\newcommand{\qpolyalg}[1][{\bsym}]{\obstuff{P}_{\qpar}^{#1}}
\newcommand{\apolyalg}[1][{\bsymv}]{\obstuff{P}_{\mathbb{A}}^{#1}}

\newcommand{\basisdet}{\setstuff{B}_{Det}}
\newcommand{\detsym}[1][{ij}]{a_{#1}}
\newcommand{\moduledet}[1][x]{\mathcal{D}^{#1}}
\newcommand{\basisgt}{\setstuff{B}_{GT}}

\newcommand{\sroot}[1][i]{\boldsymbol{\alpha}_{#1}}

\newcommand{\braid}[1][n]{\setstuff{B}_{#1}}
\newcommand{\pbraid}[1][n]{\setstuff{PB}_{#1}}
\newcommand{\ibraid}[2][n]{\setstuff{B}_{#1}^{#2}}
\newcommand{\bgen}[1][i]{\beta_{#1}}

\newcommand{\infbraid}[1][n]{\setstuff{PB}_{#1}^{\varepsilon}}
\newcommand{\infbgen}[1][{ij}]{\beta_{#1}^{\varepsilon}}
\newcommand{\infogen}[1][{ij}]{\gamma_{#1}^{\varepsilon}}

\newcommand{\casimir}[1][{i,j}]{Cas_{#1}}
\newcommand{\casimirsmall}[1][{i,j}]{cas_{#1}}
\newcommand{\casalg}[1][n]{\setstuff{Cas}_{#1}}
\newcommand{\ealg}[1][n]{\setstuff{PE}_{#1}^{\varepsilon}}
\newcommand{\evalg}[1][n]{\setstuff{PE}_{#1,\vpar}^{\varepsilon}}
\newcommand{\eqalg}[1][n]{\setstuff{PE}_{#1,\qpar}^{\varepsilon}}

\newcommand{\qlkb}[2][{\qpar,\bsym}]{\obstuff{LKB}_{#1}^{#2}}

\newcommand{\alkb}[2][{\mathbb{A},\bsymv}]{\obstuff{LKB}_{#1}^{#2}}

\newcommand{\lkbscale}[1][{+1}]{t^{#1}}

\newcommand{\rmatrix}[1][{\para_{1},\para_{2}}]{\check{r}_{#1}}
\newcommand{\rimatrix}[1][{\para_{1},\para_{2}}]{\check{r}_{#1}^{-1}}

\usepackage[all]{xy}
\usepackage{tikz}
\usetikzlibrary{cd}
\usetikzlibrary{decorations}
\usetikzlibrary{decorations.markings}
\usetikzlibrary{decorations.pathreplacing}
\usetikzlibrary{decorations.pathmorphing}
\usetikzlibrary{arrows.meta,shapes,positioning,matrix,calc}
\usetikzlibrary{shapes.callouts}

\tikzset{
	anchorbase/.style={baseline={([yshift=#1]current bounding box.center)}},
	anchorbase/.default={-0.5ex},
	tinynodes/.style={font=\tiny,text height=0.25ex,text depth=0.05ex},
	smallnodes/.style={font=\scriptsize,text height=0.75ex,text depth=0.15ex},
	mor/.style={line width=0.75,color=black,fill=cream},
	usual/.style={line width=1.2,color=black},
	crossline/.style={preaction={draw=white,line width=5.0pt,-},preaction={draw=black,line width=0.9pt,-}},
	dot/.style = {
		decoration={markings,
			post length=0.25mm,
			pre length=0.25mm,
			mark=at position #1 with {\node[circle,radius=0.15cm,inner sep=-1.2pt,color=black,fill=black]{};}
		},
		postaction={decorate}
	},
	dot/.default=1,
}
\tikzstyle directed=[postaction={decorate,decoration={markings,
		mark=at position #1 with {\arrow[line width=0.25mm, black]{>}}}}]

\usepackage{aliascnt,etoolbox}
\def\NewTheorem#1{%
	\newaliascnt{#1}{equation}%
	\newtheorem{#1}[#1]{#1}%
	\aliascntresetthe{#1}%
	\expandafter\def\csname #1autorefname\endcsname{#1}%
}
\def\equationautorefname~#1\null{(#1)\null}

\numberwithin{equation}{subsection}

\NewTheorem{Proposition}
\NewTheorem{Theorem}
\NewTheorem{Corollary}
\AtEndEnvironment{Corollary}{\null\hfill$\square$}%
\NewTheorem{Lemma}
\theoremstyle{definition}
\NewTheorem{Definition}
\NewTheorem{Notation}
\NewTheorem{Example}
\AtEndEnvironment{Example}{\null\hfill$\Diamond$}%
\theoremstyle{remark}
\NewTheorem{Remark}
\NewTheorem{Assumption}
\NewTheorem{Hypothesis}
\NewTheorem{Conjecture}
\NewTheorem{Question}
\NewTheorem{Observation}
\NewTheorem{Fact}
\NewTheorem{Conclusion}
\NewTheorem{lemma}

\usepackage[hypertexnames=false]{hyperref}
\usepackage{bookmark}
\hypersetup{
	pdftoolbar=true,
	pdfmenubar=true,
	pdffitwindow=false,
	pdfstartview={FitH},
	pdftitle={Verma Howe duality and LKB representations},
	pdfauthor={Abel Lacabanne, Daniel Tubbenhauer and Pedro Vaz},
	pdfsubject={},
	pdfcreator={Abel Lacabanne, Daniel Tubbenhauer and Pedro Vaz},
	pdfproducer={Abel Lacabanne, Daniel Tubbenhauer and Pedro Vaz},
	pdfkeywords={},
	pdfnewwindow=true,
	colorlinks=true,
	linkcolor=mydarkblue,
	citecolor=teal,
	filecolor=magenta,
	urlcolor=orchid,
	linkbordercolor=lava,
	citebordercolor=teal,
	urlbordercolor=orchid,
	linktocpage=true
}

\newcommand{\nnfootnote}[1]{%
	\begin{NoHyper}
		\renewcommand\thefootnote{}\footnote{#1}%
		\addtocounter{footnote}{-1}%
	\end{NoHyper}
}

\def\makeautorefname#1#2{\csdef{#1autorefname}{#2}}

\makeautorefname{section}{Section}%
\makeautorefname{subsection}{Section}%
\makeautorefname{subsubsection}{Section}%

\begin{document}
	\title[Verma Howe duality and LKB representations]{Verma Howe duality and LKB representations}
	\author[A. Lacabanne, D. Tubbenhauer and P. Vaz]{Abel Lacabanne, Daniel Tubbenhauer and Pedro Vaz}
	
	\address{A.L.: Laboratoire de Math{\'e}matiques Blaise Pascal (UMR 6620), Universit{\'e} Clermont Auvergne, Campus Universitaire des C{\'e}zeaux, 3 place Vasarely, 63178 Aubi{\`e}re Cedex, France,\newline \href{http://www.normalesup.org/~lacabanne}{www.normalesup.org/$\sim$lacabanne},
		\href{https://orcid.org/0000-0001-8691-3270}{ORCID 0000-0001-8691-3270}}
	\email{abel.lacabanne@uca.fr}
	
	\address{D.T.: The University of Sydney, School of Mathematics and Statistics F07, Office Carslaw 827, NSW 2006, Australia, \href{http://www.dtubbenhauer.com}{www.dtubbenhauer.com}, \href{https://orcid.org/0000-0001-7265-5047}{ORCID 0000-0001-7265-5047}}
	\email{daniel.tubbenhauer@sydney.edu.au}
	
	\address{P.V.: Institut de Recherche en Math{\'e}matique et Physique, 
		Universit{\'e} catholique de Louvain, Chemin du Cyclotron 2,  
		1348 Louvain-la-Neuve, Belgium, \href{https://perso.uclouvain.be/pedro.vaz}{https://perso.uclouvain.be/pedro.vaz}, \href{https://orcid.org/0000-0001-9422-4707}{ORCID 0000-0001-9422-4707}}
	\email{pedro.vaz@uclouvain.be}
	
	\begin{abstract}
		We establish a version of Howe duality 
		that involves a tensor product of Verma modules. Surprisingly, this duality 
			leaves the realm of lowest and highest weight modules.
		
		We quantize this duality, and as an application, we
			prove that the 
			(colored higher) LKB representations arise from this duality
			and use this description to show that they
			are simple as modules 
			for the braid group and for various of its subgroups, including 
			the pure braid group.
	\end{abstract}
	
	\nnfootnote{\textit{Mathematics Subject Classification 2020.} Primary: 17B10, 17B37; Secondary: 20C15, 20F36.}
	\nnfootnote{\textit{Keywords.} Howe duality, Verma modules, dense modules, 
		LKB representations of braid groups.}
	
	\addtocontents{toc}{\protect\setcounter{tocdepth}{1}}
	
	\maketitle
	
	\tableofcontents

%%%%%%%%%%%%%%%%%%%%%%%%%%%%%%%%%%%%%%%%%

\section{Introduction}\label{S:Introduction}

%%%%%%%%%%%%%%%%%%%%%%%%%%%%%%%%%%%%%%%%%

Arguably the most classical form of \emph{Howe duality} relates commuting 
actions of $\glgroup[m]$ and $\glgroup[n]$ on the symmetric 
algebra of $\C^{m}\otimes\C^{n}$, see \cite{Ho-remarks-invariant-theory}
or \cite{Ho-perspectives-invariant-theory}.

Howe's approach turned out to be a game changer, even in fields
beyond representation theory. For example, 
quantum versions of these dualities 
provide powerful and categorification-friendly 
descriptions of quantum invariants 
such as the colored Jones polynomial.

In this paper we prove a version 
of the Howe duality above
where symmetric powers are replaced by Verma modules. 
We call this duality \emph{Verma Howe duality}.
To the best of our knowledge, Verma Howe duality is the first example 
of a Howe duality that involves modules that are not lowest or highest weight modules.
Consequently, our proofs are very different from Howe's proofs. For example, Verma Howe duality is not a ``limit'' of symmetric Howe duality but genuinely new.

Moreover, we give an application of Verma Howe duality:
after extending Verma Howe duality to quantum groups, which is fairly straightforward, we show that 
the \emph{LKB (Lawrence--Krammer--Bigelow) representations}
and their colored and higher counterparts arise from 
quantum Verma Howe duality, which in turn enables us to 
show that the LKB representations are simple modules of various subgroups of Artin's braid group, including pure and handlebody braid groups. One direct advantage of our approach is that we can work over an arbitrary field and with a large variety of the involved parameters.

%%%%%%%%%%%%%%%%%%%%%%%%%%%%%%%%%%%%%%%%%

\subsection{Schur--Weyl(--Brauer) and Howe dualities}\label{SS:IntroA}

%%%%%%%%%%%%%%%%%%%%%%%%%%%%%%%%%%%%%%%%%

Three main themes in Weyl's seminal book 
\emph{``The classical groups''} \cite{We-classical-groups}
are the study of polynomial invariants for 
actions of the eponymous classical 
groups, and, more or less equivalent, 
decomposition of the tensor algebra for such an action, and, 
again more or less equivalent, the description of the 
invariants in the tensor algebra.

The two most prominent examples that fit into Weyl's 
setting are the celebrated \emph{Schur--Weyl duality} 
\cite{Sc-schur-weyl}
for tensor invariants of $\glgroup$ and 
\emph{Brauer duality} \cite{Br-brauer-algebra-original}
for tensor invariants of 
$\ogroup$ and 
$\spgroup$ (for the symplectic group 
$m$ is even). 
Both of these were studied by using commuting actions of $\glgroup$, and
$\ogroup$, $\spgroup$ on one side and 
the symmetric group $S_{n}$ and the Brauer algebra, respectively, 
on the other side, both acting on a tensor product 
of the defining representation of the classical groups in question. 
In this commuting-action-approach, for example
Schur--Weyl duality essentially reads:
\begin{enumerate}[label=(\Alph*)]

\item There are commuting actions of $\glgroup$ and $S_{n}$ 
on $(\C^{m})^{\otimes n}$.

\item The two actions generate each others centralizer.

\item The $\glgroup$-$S_{n}$ bimodule $(\C^{m})^{\otimes n}$ 
can be explicitly decomposed into a direct sum of nonisomorphic
simple $\glgroup$ modules tensored with nonisomorphic simple 
$S_{n}$ modules.

\end{enumerate}
A statement of this form is what we call a 
\emph{double centralizer} ({\aka} double commutant) approach.

Howe \cite{Ho-remarks-invariant-theory}, \cite{Ho-perspectives-invariant-theory} studied 
polynomial invariants, {\eg} via symmetric powers, 
of classical groups using a 
double centralizer 
approach, and the resulting dualities 
are called \emph{Howe dualities} in this paper.
A prominent example is \emph{symmetric Howe 
duality} where $\glgroup[m]$ and $\glgroup[n]$ 
act on $\medoplus_{k\in\N}\mathrm{Sym}^{k}(\C^{m}\otimes\C^{n})$. 
Howe, albeit formulated differently, proves (A)-(C) as above for this 
and other dualities.

It is not surprising that Howe-type dualities have 
been of paramount importance for the representation theory of
reductive groups ever since, see also 
\cite{ChWa-dualities-super} for a summary of various 
such dualities, but are also pervasive in other fields.
For example, in the early stages of quantum group theory 
Jimbo studied \emph{quantum Schur--Weyl duality} \cite{Ji-q-schur-weyl}, 
which, in one way or the other, is central for the study 
of quantum invariants: that the Jones polynomial 
arises from the Temperley--Lieb calculus \cite{Jo-jones-polynomial} 
is an instance of quantum Schur--Weyl duality, although 
originally not formulated as such. And this is just the tip of the iceberg.

It did not take long for \emph{quantum Howe dualities} to 
appear, see e.g. \cite{NoUmWa-sl2-son-duality} for an early reference.
Also due to their relation to diagrammatics,
quantum Howe dualities have been studied intensively since their 
first appearance in the 1990s, and also turned 
out to be very useful for the study 
of quantum invariants. For a few type $A$ examples of such quantum 
Howe dualities, 
see \cite{LeZhZh-q-first-fundamental-theorem}, 
\cite{CaKaMo-webs-skew-howe}
for quantum exterior and \cite{RoTu-symmetric-howe} 
for quantum symmetric Howe duality, and 
for some more ``exotic type $A$ settings'', 
see \cite{QuSa-mixed-skew-howe}, 
\cite{TuVaWe-super-howe} or \cite{ChWa-quantum-q-howe}, 
\cite{BrDaKu-quantum-type-q-webs}.

\begin{Remark}\label{R:IntroductionCoideals}
Quantum Howe dualities are of course not restricted to type $A$, 
but the reader should be warned at this stage: 
experience tells us that
quantum Howe dualities 
often run into quantization issues and nonstandard 
quantum objects tend to pop up. Examples are 
\cite{NoUmWa-sl2-son-duality}, \cite{EhSt-nw-algebras-howe} or \cite{SaTu-bcd-webs} where coideal subalgebras as in \cite{NoSu-q-symmetric-spaces} appear. There are even such phenomena
that are entirely in type $A$ see {\eg} \cite{LaTuVa-annular-webs-levi} 
and related quantization issues in
\cite{CaKa-q-satake-sln}, \cite{QuWe-extremal-projectors-2}.

Quantum exterior and symmetric Howe dualities as well as their Verma 
counterparts are notable exceptions, and 
the quantization in these cases is not a big deal. In fact, our 
proofs will mostly stay in the non-quantum setting and 
the quantum case then follows using a flatness argument.
\end{Remark}

%%%%%%%%%%%%%%%%%%%%%%%%%%%%%%%%%%%%%%%%%

\subsection{What this paper does}\label{SS:IntroB}

%%%%%%%%%%%%%%%%%%%%%%%%%%%%%%%%%%%%%%%%%

The main theorem of this 
paper is \autoref{T:DualityHowe} where we formulate 
a \emph{(quantum) Verma Howe duality}.
To explain the main points let us be less general than 
\autoref{T:DualityHowe} actually is. For example, as we wrote in \autoref{R:IntroductionCoideals}, 
quantization is not an issue for us 
and we can work over with general fields and quite general 
parameters, see \autoref{R:IntroductionField}, but we stay in the classical 
case in this introduction for simplicity.
The classical, non-quantum, version of \autoref{T:DualityHowe}
is then still more general than the following.

To work with Verma modules we go from the Lie group 
to the Lie algebras.
For generic enough $\para_{i}\in\C$, 
where $i\in\{1,\dots,n\}$, let $\verma{\para_{i}}$ be 
the $\ugl[2]$ Verma module of highest $\sln[2]$ weight 
$\para_{i}$. 
We take the tensor product 
$\verma{\para_{1}}\otimes\dots\otimes\verma{\para_{n}}$.
For the same reason as for symmetric Howe duality, we then take a certain 
direct sum of the $\verma{\para_{1}}\otimes\dots\otimes\verma{\para_{n}}$. Call this direct sum $\verma{\oplus\bsym}$ where $\bsym=(\para_{1},\dots,\para_{n})$.

Now, essentially by definition, $\ugl[2]$ acts 
on $\verma{\oplus\bsym}$ and we 
also construct a dual action of $\ugl[n]$ 
on $\verma{\oplus\bsym}$. 
Using the double centralizer approach, \autoref{T:DualityHowe} 
states and proves (A)-(C) for the $\ugl[2]$-$\ugl[n]$ 
bimodule $\verma{\oplus\bsym}$.

Since all symmetric powers for $\ugl[2]$ are quotients of Verma modules,
we think of this Verma Howe duality as a generalization of 
symmetric Howe duality (with a caveat, see \autoref{SS:IntroC} below).
Verma Howe duality is however much more difficult to prove: 
Firstly, the whole setting is, by its very nature, infinite dimensional 
and most of the classical 
statements need to be appropriately reformulated and adjusted 
to the infinite dimensional setting.
Second, and more importantly, the simple $\ugl[n]$ appearing in (C) 
are neither highest nor lowest weight modules; they are 
simple \emph{dense} (weight) modules in the sense of \cite{Ma-classification-weight-modules}.
This is, to the best of our knowledge, very different from all other Howe-type dualities in the literature and makes calculations (for example actions of Casimir 
elements) much more involved. 
In particular, we need to take quite a detour to identify the dense modules 
explicitly and we crucially use results from \cite{Ma-gt-dense} 
and \cite{MiToLa-casimir-braids}, and implicitly computer help, 
to identify them.
(This is also our main reason to stay 
with $\ugl[2]$ instead of $\ugl[m]$.) Along the way 
we partially generalize \cite{Ma-gt-dense} so that we can use fairly general parameters.

\begin{Remark}\label{R:IntroductionField}
Let us also stress that our approach works in quite some generality. That is, we work over an arbitrary field $\K$ and fix a quantum parameter that is not a root of unity. Moreover, the $\para_{i}$ of the Verma modules are, up to a certain degree, allowed to be integers, see \autoref{D:DualityGeneric} for a precise condition.
\end{Remark}

As an application of \autoref{T:DualityHowe} we prove that 
the (colored higher) LKB representations constructed in
\cite{JaKe-verma-lkb} and \cite{Ma-colored-lkb} are simple 
as modules of the associated (colored) braid groups. This 
not just gives a new proof of 
\cite[Theorem 3]{JaKe-verma-lkb} but also strengthens
the result of Jackson--Kerler quite a bit: we prove simplicity for much smaller groups, namely the corresponding pure braid groups. Moreover, Jackson--Kerler work over $\K(q)$ for $\Q\subset\K$ 
and with a generic 
parameter for the LKB representations. Our setting is more general, see \autoref{R:IntroductionField}. In fact, we think it is remarkable that the LKB representations stay simple even after specializing some parameters or leaving characteristic zero. 
Finally, since we can allow different parameters, our methods also relate the LKB representations to handlebody braid groups as in {\eg} 
\cite{Ve-handlebodies}, \cite{HaOlLa-handlebodies}, \cite{RoTu-homflypt-handlebody} or \cite{TuVa-handlebody}.

%%%%%%%%%%%%%%%%%%%%%%%%%%%%%%%%%%%%%%%%%

\subsection{Outlook}\label{SS:IntroC}

%%%%%%%%%%%%%%%%%%%%%%%%%%%%%%%%%%%%%%%%%

Separate from the evident question how to replace $\ugl[2]$ 
by $\ugl[m]$, here are a few directions one could try to explore:
\begin{enumerate}

\item While (A) and (B) 
as above often hold in more generality, (C) is using 
that the underlying representation is semisimple.
The nonsemisimple versions of some of the above are known, 
see for example \cite{DuPaSc-schur-weyl-tilting} 
for an integral version of 
quantum Schur--Weyl duality. But these are
also much more involved and often need 
some form of tilting theory.

A nonsemisimple version \autoref{T:DualityHowe} would be a true 
generalization of quantum symmetric Howe duality since the cases 
where only symmetric powers appear within the Verma modules are precisely ruled out by our condition in \autoref{D:DualityGeneric}. However, we 
can still have symmetric powers but need at least also a ``generic enough'' 
highest weight.

\item Several 
papers discuss dualities involving one Verma and tensor 
products of finite dimensional 
modules, see {\eg} \cite{IoLeZh-verma-schur-weyl} or \cite{LaVa-schur-weyl-ariki-koike}.
It would be interesting to compare these to this work, also with 
an eye on categorification of the story as in \cite{LaNaVa-tensor-product-blob}.

\item Another interesting direction is the identification of 
the LKB representations with specialized parameters as cell representations 
of algebras within the symmetric web category from \cite{RoTu-symmetric-howe}.
We suspect that this is a consequence of Verma Howe duality for 
(the quantum version of) $\para_{i}\in\N$.
Note that special cases of this are known: Jones' work 
\cite{Jo-jones-polynomial} implicitly showed
the respective statement for the Burau representation and the Temperley--Lieb 
calculus, and \cite{Zi-lkb-bmw} implicitly showed an analog for the LKB representation.
Note that Temperley--Lieb and the Brauer-type calculus used 
in \cite{Jo-jones-polynomial} and 
\cite{Zi-lkb-bmw}, respectively, are special cases of the 
symmetric web calculus. (For the Temperley--Lieb calculus 
this is clear, while 
the Brauer-type calculus makes its appearance due to the 
``small number coincidence'' that matches $\sogroup$ representations and odd dimensional $\slgroup[2]$ representations.)
In \cite[Lemma 6]{Fo-braid-group-reps} it is shown that the reduced Burau representation of the $n$ strand braid group is simple if and only if quantum $n$ does not vanish, and Verma Howe duality should be helpful to prove similar results for the other LKB representations.

\item A striking question is how to categorify Verma Howe duality. 
We suspect this should be related to categorification of tensor products 
of infinite dimensional representations as in 
\cite{DuNa-categorical-lkb}. One could also hope, in some sense, that a categorification of 
LKB representations would be an upshot of such a categorical Verma Howe duality.

\end{enumerate}
\medskip

\noindent\textbf{Acknowledgments.}
We like to thank Volodymyr Mazorchuk for many helpful exchanges of emails, 
and for explaining various properties of dense modules to us. We also thank the referee for helpful comments.
Part of this paper were done after having consulted
Magma and Mathematica. Their help is gratefully acknowledged.

%%%%%%%%%%%%%%%%%%%%%%%%%%%%%%%%%%%%%%%%%

\section{A duality involving Verma modules}\label{S:Duality}

%%%%%%%%%%%%%%%%%%%%%%%%%%%%%%%%%%%%%%%%%

In this section we state a duality that we call \emph{(quantum) Verma Howe duality}.

\begin{Remark}\label{R:DualityColors}
We use colors in this paper, but these are a visual aid 
and do not have other significance. In particular, the paper is readable 
in black-and-white without restrictions.
\end{Remark}

%%%%%%%%%%%%%%%%%%%%%%%%%%%%%%%%%%%%%%%%%

\subsection{Verma and dense modules}\label{SS:DualityVermaDense}

%%%%%%%%%%%%%%%%%%%%%%%%%%%%%%%%%%%%%%%%%

The following specifies the underlying field:

\begin{Notation}\label{N:DualityField}
Fix an arbitrary field $\K$
and an element $\qpar\in\K\setminus\{0\}$ that is not of finite order. We call $\qpar$ the \emph{quantum parameter}.

We additionally allow $q=\pm 1$, but then we assume that $\K$ is of characteristic zero. This is the \emph{non-quantum} or \emph{classical case}. The reader is warned that the below is tailor-made for the quantum case and needs to be adjusted for the classical case. We leave 
the adjustments to the reader.
\end{Notation}

We consider the \emph{quantum enveloping algebra $\uqgl[2]$ of $\gl[2]$}
over $\K$ with respect to the quantum parameter $\qpar$. We 
specify our conventions later on in \autoref{SS:ProofQuantum} and for now it is enough to 
know that $\uqgl[2]$ is, as a $\K$ algebra, 
generated by $E$, $F$, $L_{1}^{\pm 1}$ and 
$L_{2}^{\pm 1}$.

From now on fix $n\in\Z_{\geq 1}$.

\begin{Notation}\label{N:DualityBold}
\leavevmode
\begin{enumerate}

\item We use a bold font for tuples, 
{\eg} $\bsym=(\para_{1},\dots,\para_{n})\in\K^{n}$.

\item Whenever an index of tuple is not defined 
but appears in a formula, then 
the associated element is zero, by convention. 
For example, $\para_{<1}=\para_{>n}=0$ if we specify 
$\bsym=(\para_{1},\dots,\para_{n})$.

\item We will also use sums of the form 
$a_{1}+a_{2}+\dots+a_{k-1}+a_{k}$ for $k\in\Z_{\geq 1}$ 
very often in this paper, 
and we abbreviate them to $\sums{a_{k}}=\medsum_{i=1}^{k}a_{i}$.

\item Denote by $\stdvec=(0,\dots,0,1,0,\dots,0)$ the 
tuple with the $i$th entry being $1$, and $\sroot[i]=\stdvec-\stdvec[i+1]$.
We also use $\stdvec[{ij}]$ below meaning a matrix-style notation with only 
one nonzero entry.

\end{enumerate}
\end{Notation}

\begin{Definition}\label{D:DualityQNum}
Given $\para\in\K$ we consider the field $\Kpara[{\para}]=\K(\qpar^{\para})$.
We define the \emph{quantum numbers} as
$\qnum{x}=\frac{\qpar^{x}-\qpar^{-x}}
{\qpar-\qpar^{-1}}\in\Kpara[{\para}]$,
where $x\in\Z$ or $x\in\para+\Z$. Similarly, for
$\bsym\in\K^{n}$ we use the field 
$\Kpara=\K(\qpar^{\para_{1}},\dots,\qpar^{\para_{n}})$ and quantum numbers 
will be elements of $\Kpara$.
\end{Definition}

The following will be often used silently throughout:

\begin{Lemma}\label{L:DualityField}
All quantum numbers are nonzero and thus invertible.
\end{Lemma}

\begin{proof}
Easy since \autoref{N:DualityField} forces this to be true, in particular, we need characteristic zero for $q=\pm 1$. Details are omitted.
\end{proof}

\begin{Remark}\label{R:DualityParameters}
The tuple $\bsym\in\K^{n}$ consist of the 
underlying \emph{parameters} that we use. 
Note that our formulation includes the case where 
the quantum parameter $\qpar$ 
and the $\para_{i}$ are formal variables by {\eg} choosing $\K=\Q(Z,Z_{1},\dots,Z_{n})$, 
for indeterminates $Z$ and $Z_{i}$, and $\qpar=Z$, $\para_{i}=Z_{i}$.
In contrast, the parameters could be in $\Z\subset\K$, but we 
partially need to avoid that, see {\eg} \autoref{D:DualityGeneric} below.
It is allowed that some (or even all) 
of the $\para_{i}$ are the same.
\end{Remark}

We consider $\uqgl[2]$ also over 
fields such as $\Kpara$ by scalar extension.
The parameters only play a role for 
$\uqgl[2]$ modules and not for $\uqgl[2]$ itself.

\begin{Definition}\label{D:DualityVerma}
For any $\para\in\K$ the \emph{(quantum) dual Verma module} $\qverma{\para}$ 
of highest weight $\lambda$ is $\qverma{\para}=\Kpara[{\para}]\{m_{i}|i\in\N\}$ 
as a $\Kpara[{\para}]$ vector space and the 
left $\uqgl[2]$ action is
\begin{gather}\label{Eq:DualityAction1}
\begin{gathered}
E\acts m_{i}=\qnum{i}\cdot m_{i-1}
,\quad
F\acts m_{i}=\qnum{\para-i}\cdot m_{i+1}
,\\
L_{1}\acts m_{i}=\qpar^{\para-i}\cdot m_{i}
,\quad
L_{2}\acts m_{i}=\qpar^{i}\cdot m_{i}
,
\end{gathered}
\end{gather}
where we use the quantum numbers from \autoref{D:DualityQNum} 
and let $m_{-1}=0$.

More generally, we define $\qverma{{\para,t}}$ for $t\in\K$ 
by tensoring $\qverma{\para-2t}$ with the one dimensional $\uqgl[2]$ module of 
highest $\gl[2]$ weight $(t,t)$.
\end{Definition}

We call the $\qverma{\para}$ 
\emph{Verma modules} for simplicity although 
they coincide with what are often called dual Verma modules in the literature,
{\eg} our modules correspond to $M^{\vee}$ in \cite{Hu-cat-o}.

\begin{Remark}\label{R:DualityVerma}
The highest weight of $\qverma{\para}$ is strictly 
speaking $(\qpar^{\para},\qpar^{0})$, but in 
\autoref{D:DualityVerma}, and throughout, we use
the notion \emph{weight} in
the sense of classical $\sln[2]$ weight combinatorics.
We however sometimes need to be more specific. For example,
for $\qverma{{\para,t}}$ we need 
the $\gl[2]$ weight notation and the classical highest weight 
of $\qverma{{\para,t}}$ is $(\para-t,t)$ which is the same 
as $(\para-2t,0)$ when restricted to $\sln[2]$ weight notation. 
Whenever we use $\gl[2]$ notation we point that out.
\end{Remark}

\begin{Example}\label{E:DualityVerma}
The $\uqgl[2]$ module $\qverma{\para}$ is given by the usual picture but with 
slightly more generic quantum numbers in the action:
\begin{gather*}
\scalebox{0.8}{$\qverma{\para}\leftrightsquigarrow
\begin{tikzpicture}[anchorbase]
\foreach \x in {0,...,5} {
\node at (-2*\x,0){$m_{\x}$};
\draw[spinach,->] (-2*\x-0.25,0.2) to[out=135,in=45] node[above]{$\scalebox{0.8}{$\qnum{\para{-}\x}$}$} (-2*\x-1.75,0.2);
\draw[tomato,->] (-2*\x-1.75+2,-0.2) to[out=315,in=225] node[below]{\scalebox{0.8}{$\qnum{\x}$}} (-2*\x-0.25+2,-0.2);
\draw[orchid,->] (-2*\x-0.1,0.2) to[out=135,in=180] (-2*\x,0.7)node[above]{$\qpar^{\para-2\cdot\x}$} to[out=0,in=45] (-2*\x+0.1,0.2);
}
\draw[tomato,->] (-2*6-1.75+2,-0.2) to[out=315,in=225] node[below]{\scalebox{0.8}{$\qnum{6}$}} (-2*6-0.25+2,-0.2);
\node at (-2*6,0){$\cdots$};
\end{tikzpicture},$}
\\
{\color{tomato}E}\text{ moves to the right}
,\quad
{\color{spinach}F}\text{ moves to the left}
,\quad
{\color{orchid}K}=L_{1}L_{2}^{-1}\text{ is a loop}.
\end{gather*}
The highest weight is $\para$.
\end{Example}

We will also use 
the \emph{quantum enveloping algebra $\uqgl[n]$ of $\gl[n]$}, again over $\Kpara$, with 
conventions specified later on, see \autoref{SS:ProofQuantum}. The generators are 
$E_{i}$, $F_{i}$ for $i\in\set{1,\dots,n-1}$,
and $L_{i}^{\pm 1}$ for $i\in\set{1,\dots,n}$.

\begin{Notation}\label{N:DualityLeftRight}
We consider $\uqgl[2]$ and $\uqgl[n]$ as different algebras, even if $n=2$.
All $\uqgl[2]$ modules used in this paper are left
$\uqgl[2]$ modules, while all $\uqgl[n]$ modules used 
in this paper are right 
$\uqgl[n]$ modules. If we mean either $\uqgl[2]$ or $\uqgl[n]$, then 
we will write $\uqgl[k]$. 
Similarly for their classical versions, 
and we will often drop the adjectives left and right.
\end{Notation}

We will need 
certain $\uqgl[n]$ modules with bases 
indexed by:

\begin{Definition}\label{D:DualityGT}
Fix $\bsym[m]=(m_{1},\dots,m_{n})\in\K^{n}$ and $\bsym[y]=(y_{1},\dots,y_{n-1})\in\K^{n-1}$ 
such that $m_{2}-m_{3},\dots,m_{n-1}-m_{n}\in\N$ and $y_{i}-m_{2}\notin\Z$ 
for all $i\in\{1,\dots,n-1\}$. A
\emph{GT (Gelfand--Tsetlin) pattern} $\gt[\vec{x}]$
for $(\bsym[m],\bsym[y])$ is a triangular array
of the form
\begin{gather}\label{Eq:DualityGT}
\gt[\vec{x}]=
\begin{array}{CCCCCCCCC}
\cellcolor{tomato!50}$x_{n1}$ & & \cellcolor{spinach!50}$x_{n2}$ & \cellcolor{spinach!50} & \cellcolor{spinach!50}$\dots$ & \cellcolor{spinach!50} & \cellcolor{spinach!50}$\dots$ & \cellcolor{spinach!50} & \cellcolor{spinach!50}$x_{nn}$ \\
& \cellcolor{tomato!50}$\ddots$ & & \cellcolor{spinach!50}$\ddots$ & \cellcolor{spinach!50} & \cellcolor{spinach!50}$\rddots$ & \cellcolor{spinach!50} & \cellcolor{spinach!50}$\rddots$ & \\
& & \cellcolor{tomato!50}$x_{31}$ & & \cellcolor{spinach!50}$x_{32}$ & \cellcolor{spinach!50} & \cellcolor{spinach!50}$x_{33}$ & & \\
& & & \cellcolor{tomato!50}$x_{21}$ & & \cellcolor{spinach!50}$x_{22}$ & & & \\
& & & & \cellcolor{tomato!50}$x_{11}$ & & & &
\end{array}
\end{gather}
where $\vec{x}=(x_{n1},\dots,x_{nn},x_{(n-1)1},\dots)$ ({\ie} the pattern read row-wise)
is such that: 
\begin{enumerate}[label=(\roman*)]

\item $x_{ni}=m_{i}$ for $i\in\{1,\dots,n\}$ ($\bsym[m]$ gives the top row),

\item $x_{i1}-y_{i}\in\Z$ for $i\in\{1,\dots,n\}$ ($\bsym[y]$ gives the first diagonal up to integers),

\item $x_{jk}-x_{(j-1)k}\in\N$ and $x_{(j-1)k}-x_{j(k+1)}\in\N$ for 
$j\in\{3,\dots,n\}$ and $k\in\{2,\dots,j\}$, that is,
\begin{gather*}
\begin{array}{CCCCC}
$x_{jk}$ & & $\geq$ & & $x_{j(k+1)}$ \\
& \raisebox{0.05cm}{\rotatebox{-45}{$\geq$}} & & \raisebox{-0.2cm}{\rotatebox{45}{$\geq$}} & \\
& & $x_{(j-1)k}$ & &
\end{array}
.
\end{gather*}

\end{enumerate}

A \emph{two diagonal GT pattern} (appearing in Verma Howe duality)
is a GT pattern with $\bsym[m]=(x_{n}=\sums{\para_{n}}+b,c,0,\dots,0)$ 
with $b,c\in\Z$ to be chosen and $\bsym[y]$ determined by $y_{i}=\sums{\para_{i}}$
for $i\in\{1,\dots,n-1\}$.
Letting $c_{n}=c$, we denote these by
\begin{gather}\label{Eq:DualityTwoDiagonal}
\gt[{\bsym[x],\bsym[c]}]=
\begin{array}{CCCCCCCCC}
\cellcolor{tomato!50}$x_{n}$ & & \cellcolor{spinach!50}$c_{n}$ & \cellcolor{spinach!50} & \cellcolor{spinach!50}$0$ & \cellcolor{spinach!50} & \cellcolor{spinach!50}$\dots$ & \cellcolor{spinach!50} & \cellcolor{spinach!50}$0$ \\
& \cellcolor{tomato!50}$x_{n{-}1}$ & & \cellcolor{spinach!50}$c_{n{-}1}$ & \cellcolor{spinach!50} & \cellcolor{spinach!50}$\ddots$ & \cellcolor{spinach!50} & \cellcolor{spinach!50}$\rddots$ & \\
& & \cellcolor{tomato!50}$\ddots$ & & \cellcolor{spinach!50}$\ddots$ & \cellcolor{spinach!50} & \cellcolor{spinach!50}$0$ & & \\
& & & \cellcolor{tomato!50}$x_{2}$ & & \cellcolor{spinach!50}$c_{2}$ & & & \\
& & & & \cellcolor{tomato!50}$x_{1}$ & & & &
\end{array}
,
\end{gather}
with $\bsym[x]=(x_{1},\dots,x_{n})\in(\Kpara[{\para}])^{n}$ 
as above and $\bsym[c]=(c_{2},\dots,c_{n})\in\N^{n-1}$.
\end{Definition}

In \autoref{Eq:DualityGT} and \autoref{Eq:DualityTwoDiagonal} 
we have shaded the parts of 
the GT patterns which play significantly different roles. 
We call the shaded block to the right 
the \emph{integral part} of the pattern 
since all the patterns we need will have integral entries in 
this part, and nonintegral entries otherwise.

\begin{Definition}\label{D:DualityGTModule}
Consider two diagonal GT patterns. 
%as in \autoref{D:DualityGT}.
Define the \emph{dense module}, as a $\Kpara[{\para}]$ vector space, 
as \[\qdense{\bsym[m],\bsym[y]}=\Kpara[{\para}]
\{\gt[\vec{x}]|\gt[\vec{x}]\text{ is a two diagonal GT pattern for }(\bsym[m],\bsym[y])\}\] and the 
$\uqgl[n]$ action is given later in \autoref{Eq:ProofDiffOpQuantumTwo} 
(on a different basis).
\end{Definition}

\begin{Example}\label{E:DualityGTModule}
If $n=2$, then the only entry in a GT patterns 
that is not completely determined by $(\bsym[m],\bsym[y])$ 
is $x_{11}$. The latter is some integer shift of $y_{1}$, so we 
can index a basis of $\qdense{\bsym[m],\bsym[y]}$ as $\{w_{i}|i\in\Z\}$.
For certain values of $A_{i}$, 
$B_{i}$ and $C_{i}$ that can be explicitly obtained from 
\autoref{Eq:ProofDiffOpQuantumTwo} the picture is then
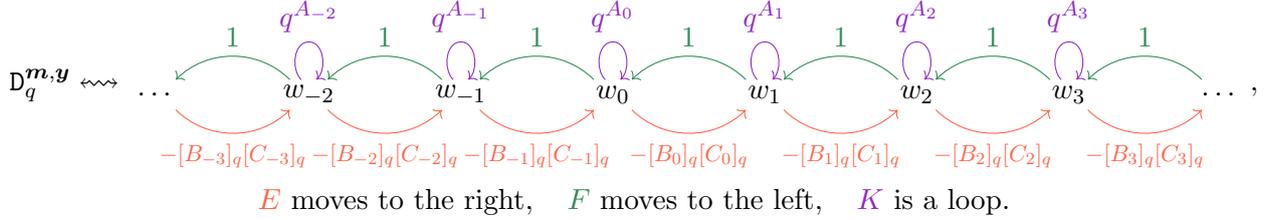
\begin{gather*}
\scalebox{0.75}{$\qdense{\bsym[m],\bsym[y]}\leftrightsquigarrow
\begin{tikzpicture}[anchorbase]
\foreach \x [evaluate=\x as \xx using int(-\x+3)] in {0,...,5} {
\node at (-2*\x,0){$w_{\xx}$};
\draw[spinach,->] (-2*\x-0.25,0.2) to[out=135,in=45] node[above]{$1$} (-2*\x-1.75,0.2);
\draw[tomato,->] (-2*\x-1.75+2,-0.2) to[out=315,in=225] node[below]{\scalebox{0.8}{$-\qnum{B_{\xx}}\qnum{C_{\xx}}$}} (-2*\x-0.25+2,-0.2);
\draw[orchid,->] (-2*\x-0.1,0.2) to[out=135,in=180] (-2*\x,0.7)node[above]{$\qpar^{A_{\xx}}$} to[out=0,in=45] (-2*\x+0.1,0.2);
}
\draw[spinach,->] (-2*0-0.25+2,0.2) to[out=135,in=45] node[above]{$1$} (-2*0-1.75+2,0.2);
\draw[tomato,->] (-2*6-1.75+2,-0.2) to[out=315,in=225] node[below]{\scalebox{0.8}{$-\qnum{B_{-3}}\qnum{C_{-3}}$}} (-2*6-0.25+2,-0.2);
\node at (-2*6,0){$\cdots$};
\node at (2*1,0){$\cdots$};
\end{tikzpicture},$}
\\
{\color{tomato}E}\text{ moves to the right}
,\quad
{\color{spinach}F}\text{ moves to the left}
,\quad
{\color{orchid}K}\text{ is a loop}.
\end{gather*}
The module $\qdense{\bsym[m],\bsym[y]}$ has neither a highest nor a lowest weight.
Note also that the conventions for the scalars in this example are different 
from \autoref{E:DualityVerma}. (That is also why the basis vectors here are denoted by $w_{i}$ and not by $m_{i}$.) 
But that can be fixed by appropriate base change. 
\end{Example}

\begin{Remark}\label{R:DualityGTModule}
For $\uqgl[2]$ there are four \emph{interval-type pictures} 
as in \autoref{E:DualityVerma} and \autoref{E:DualityGTModule}. 
First, a finite interval $[a,b]$, having a highest 
and a lowest weight, which corresponds to a finite dimensional 
$\uqgl[2]$ module. One could also use $]{-}\infty,b]$ 
or $[a,\infty[$, and the associated $\uqgl[2]$ modules are 
Verma and coVerma modules, respectively. These have either a 
highest or a lowest weight. Finally, the interval 
$]{-}\infty,\infty[=\R$ corresponds to the dense modules 
and these have neither a highest 
nor a lowest weight. In this sense, dense modules 
are a natural family of $\uqgl[2]$ modules.

More generally, dense modules appear in the study of 
weight modules for $\uqgl[k]$. That is, every simple 
weight module of $\uqgl[k]$ is dense or induced from a dense 
module, see \cite{Fu-weight-modules} and \cite{Fe-weight-modules},
which reduces the classification of simple 
weight modules to dense modules. 
Hence, one could say that 
dense modules are prototypical weight modules.
\end{Remark}

\begin{Lemma}\label{L:DualityAction}
\autoref{Eq:DualityAction1} and 
\autoref{D:DualityGTModule}
endow $\qverma{\para}$ and 
$\qdense{\bsym[m],\bsym[y]}$, respectively, with structures 
of $\uqgl[2]$ and $\uqgl$ modules.
\end{Lemma}

\begin{proof}
Well-known and easy for $\qverma{\para}$, and this follows from 
\autoref{Eq:ProofDiffOpQuantumTwo} below for $\qdense{\bsym[m],\bsym[y]}$.
\end{proof}

\begin{Definition}\label{D:DualityGeneric}
We call $\bsym$ \emph{admissible parameters} if exists a permutation 
$\sigma\in\Aut\{1,\dots,n\}$ such that
$\sums{\para_{\sigma(k)}}\notin\Z$ for all $k\in\{1,\dots,n\}$.
\end{Definition}

\begin{Example}\label{E:DualityGeneric}
Note that \autoref{D:DualityGeneric} allows to have some $\para_{i}\in\Z$. For example, the parameters 
$\bsym=(1,2,3,\pi,4,5,6)\in\R^{7}$ are admissible.
\end{Example}

We will need admissible parameters because of 
\autoref{R:DualitySymmetric} below and also because of:

\begin{Lemma}\label{L:ProofDualityDense}
For admissible parameters we have that the $\uqgl[2]$ module $\qverma{\para}$ 
and the $\uqgl[n]$ module $\qdense{b,c}$
are simple. Similarly, $\qverma{\para,t}$ is a simple $\uqgl[2]$ module 
if $\para-2t$ is generic.
\end{Lemma}

\begin{proof}
For the dense modules 
we will show this later in \autoref{L:ProofIsSimpleTwo} while 
$\para\notin\Z$ 
implies simplicity of $\qverma{\para}$, 
as usual in the theory, {\cf} \cite[Section 1.5]{Hu-cat-o}.
\end{proof}

%%%%%%%%%%%%%%%%%%%%%%%%%%%%%%%%%%%%%%%%%

\subsection{Verma Howe duality}\label{SS:DualityHowe}

%%%%%%%%%%%%%%%%%%%%%%%%%%%%%%%%%%%%%%%%%

Since we work with infinite dimensional 
$\Kpara$ vector spaces and their homomorphisms, 
we need to be careful with respect to 
finite {\versus} infinite sums. To avoid convergence 
issues, we use the following definition, where rings, as throughout, 
are associative and unital.

\begin{Definition}\label{D:DualityInfSums}
Let $\ring[S]\subset\ring[T]$ be two 
rings, and let $\module[M]$ 
be a left (or right) $\ring[T]$ module. We call $\ring[S]$ a \emph{dense subring} of $\ring[T]$ (with respect to $\module[M]$)
if for any $t\in\ring[T]$ and $m_{1},\dots,m_{k}\in\module[M]$ there 
exists $s\in\ring[S]$ such that $s\acts m_{i}=t\acts m_{i}$ 
(or $m_{i}\acts s=m_{i}\acts t$)
for $i\in\{1,\dots,k\}$. 

We say 
$\{s_{i}|i\in I\}\subset\ring[T]$ \emph{densely-generates} $\ring[T]$ 
(with respect to a fixed $\module[M]$)
if $\{s_{i}|i\in I\}$ generates a dense subring of $\ring[T]$ 
and we write $\{s_{i}|i\in I\}\dgen\ring[T]$ in this case.
\end{Definition}

\begin{Notation}\label{N:DualityOpposite}
We also write $\End_{\ring[S]}(\module[M])$ instead of $\End_{\ring[S]^{op}}(\module[M])$, {\ie} we suppress the necessary but not enlightening appearance of the opposite ring.
\end{Notation}

We will write $\qdense{b,c}$ for the dense modules with $b,c$ as in \autoref{D:ProofGTBasis}.

Let $\ring[S]$ be a ring.
For a left or right $\ring[S]$ module $\module[M]$, 
the ring $\ring[S]^{\prime}=\End_{\ring[S]}(\module[M])$ 
is called the \emph{centralizer} of $\ring[S]$ (on $\module[M]$).
We call the following theorem 
\emph{(quantum) Verma Howe duality}:

\begin{Theorem}\label{T:DualityHowe}
\leavevmode	

\begin{enumerate}

\item There are commuting actions
\begin{gather*}
\uqgl[2]\actsleft
\qverma{\oplus\bsym}=\medoplus_{\bsym[d]\in\Z^{n}}\qverma{\para_{1}+d_{1}}
\otimes\dots\otimes\qverma{\para_{n}+d_{n}}
\actsright
\uqgl[n].
\end{gather*}

\item Let $\amap{\qpar}{k}$ be the algebra homomorphism 
induced by the $\uqgl[k]$ actions from (a). Then, for admissible parameters $\bsym$:
\begin{gather*}
\amap{\qpar}{2}\colon
\uqgl[2]\dgen\End_{\uqgl[n]}(\qverma{\oplus\bsym})
,\quad
\amap{\qpar}{n}\colon
\uqgl[n]\dgen\End_{\uqgl[2]}(\qverma{\oplus\bsym})
.
\end{gather*}
That is, the two actions densely-generate the others centralizer.	

\item For\! admissible parameters\! $\bsym$ we have the decomposition of the 
$\uqgl[2]$-$\uqgl[n]$ bimodule $\qverma{\oplus\bsym}$ into

\begin{gather}\label{Eq:DualityDecom}
\qverma{\oplus\bsym}\cong
\medoplus_{\substack{\countingpara\in\Z\\ t\in\N}}
\qverma{\sums{\para_{n}}+\countingpara-t,t}
\otimes\qdense{\countingpara-t,t}
.
\end{gather}
The various $\qverma{\sums{\para_{n}}+\countingpara-t,t}$ and $\qdense{\countingpara-t,t}$ 
are nonisomorphic simple $\uqgl[2]$ modules respectively $\uqgl[n]$ modules.

\end{enumerate}
There is also a similar statement in the non-quantum case which the 
reader can spell out easily themselves by removing all $\qpar$ above.
\end{Theorem}

The proof of \autoref{T:DualityHowe} is 
nontrivial and given in its own section, see \autoref{S:Proof} below.

\begin{Remark}\label{R:DualityInfSums}
If $\module[M]$ in \autoref{D:DualityInfSums} 
is finitely generated, then 
densely-generating the centralizer is the same 
as generating the centralizer. In this case 
\autoref{T:DualityHowe} is a classical Schur--Weyl(--Brauer) 
or Howe duality as in the introduction. The formulation above is copied from \cite[Section 3]{AnStTu-semisimple-tilting}, which also gives an
overview of Schur--Weyl(--Brauer) dualities.
\end{Remark}

\begin{Remark}\label{R:DualitySymmetric}
We suspect that \autoref{T:DualityHowe}.(b) works without  
assuming that we have admissible parameters, and we would expect 
tilting theory as in the proofs of \autoref{L:ProofSemisimple} 
and \autoref{L:ProofFlat}
below to play a major role. However, note that $\qverma{\para}$ for 
$\para\in\N$ is not tilting which makes the nonsemisimple 
situation much more delicate.
Note that
\autoref{T:DualityHowe} for $\bsym\in\N^{n}$ could be used 
to generalize \emph{(quantum) symmetric Howe duality} as 
in, for example, 
\cite[Theorem 2.1.2]{Ho-perspectives-invariant-theory} and \cite[Theorem 2.6]{RoTu-symmetric-howe}.
\end{Remark}

\begin{Remark}\label{R:DualityHowe}
The GT patterns in 
\autoref{T:DualityHowe} always have many zeros, exactly 
as in \autoref{Eq:DualityTwoDiagonal}. This is because we consider $\uqgl[2]$
and not $\uqgl[m]$ for general $m\in\Z_{\geq 1}$.
\end{Remark}

\begin{Remark}\label{Remark:DualityGlSl}
Verma Howe duality as in \autoref{T:DualityHowe} is formulated 
for $\big(\uqgl[2],\uqgl[n]\big)$. If the reader likes to work with 
the special linear group instead of the general linear group, 
then they can replace 
$\big(\uqgl[2],\uqgl[n]\big)$ with $\big(\uqsl[2],\uqgl[n]\big)$ 
or $\big(\uqgl[2],\uqsl[n]\big)$ in \autoref{T:DualityHowe}.
\end{Remark}

%%%%%%%%%%%%%%%%%%%%%%%%%%%%%%%%%%%%%%%%%

\section{The proof of Verma Howe duality}\label{S:Proof}

%%%%%%%%%%%%%%%%%%%%%%%%%%%%%%%%%%%%%%%%%

We first prove the classical version of \autoref{T:DualityHowe}, 
and then use a flatness argument to get the quantum version. Recall that 
in the classical case we assume that $\K$ is of characteristic zero.

%%%%%%%%%%%%%%%%%%%%%%%%%%%%%%%%%%%%%%%%%

\subsection{The classical case}\label{SS:ProofClassical}

%%%%%%%%%%%%%%%%%%%%%%%%%%%%%%%%%%%%%%%%%

We will need the Lie algebra $\gl[k]$ and it elements 
of the form $E_{ij}$. These are the $k{\times}k$ matrices with a
one in the $i$th row and $j$th column and zeros otherwise.

\begin{Lemma}\label{L:ProofFixOrder}
If \autoref{T:DualityHowe} holds for $\bsym\in\K^{n}$, then it holds for any permutation of $\bsym$ as well.
\end{Lemma}

\begin{proof}
This follows since the category of $\gl[k]$-representations is symmetric.
\end{proof}

\begin{Notation}\label{N:ProofSL2}
\leavevmode
\begin{enumerate}

\item We also write 
$E_{i}=E_{i(i+1)}$, $F_{i}=E_{(i+1)i}$ and 
$L_{i}=E_{ii}$.
For $\gl[2]$, 
we simplify this notation and use $E=E_{1}=\begin{psmallmatrix}0&1\\0&0\end{psmallmatrix}$ 
and $F=F_{1}=\begin{psmallmatrix}0&0\\1&0\end{psmallmatrix}$, and 
we also have $L_{1}=\begin{psmallmatrix}1&0\\0&0\end{psmallmatrix}$ 
and $L_{2}=\begin{psmallmatrix}0&0\\0&1\end{psmallmatrix}$.

\item We denote the operators used in actions by {\eg} $e_{ij}$ 
to distinguish then from the elements of the Lie algebras. 
The operators are always elements of some endomorphism space. 
The appearing operators will always be denoted using lowercase letters.

\item By \autoref{L:ProofFixOrder} we can and will assume that 
$\sums\para_{k}\notin\Z$ for all $k\in\{1,\dots,n\}$ instead of $\sums{\para_{\sigma(k)}}\notin\Z$ for all $k\in\{1,\dots,n\}$. This 
will be of importance in some of our formulas.

\end{enumerate}
\end{Notation}

We need the following realization of $\verma{\oplus\bsym}$.
Let $\K[\bsym[X]^{\pm 1},\bsym[Y]]$ be the algebra 
generated by indeterminates $\bsym[X]^{\pm 1}
=(X_{1}^{\pm 1},\dots,X_{n}^{\pm 1})$ 
and $\bsym[Y]=(Y_{1},\dots,Y_{n})$. 
We shift the exponents of the 
$X$ in $\K[\bsym[X]^{\pm 1},\bsym[Y]]$ by $\bsym$ 
so that powers of the variables 
$\bsym[X]$ and $\bsym[Y]$ are now in $\bsym+\Z^{n}$ 
and $\N^{n}$, respectively. The resulting $\K$ vector space 
is denoted by $\polyalg=\K[\bsym[X]^{\bsym+\Z^{n}},\bsym[Y]]$. 
We view $\polyalg$ as a $\K[\bsym[X]^{\pm 1},\bsym[Y]]$ bimodule, meaning 
that we allow multiplication by $X_{i}^{\pm 1}$ 
and by $Y_{i}$.
We also use $\polyalg[\para]$ defined 
similarly.

\begin{Definition}\label{D:ProofOperators}
For $i\in\{1,\dots,n\}$ we let operators $\partial_{X_{i}}$ and 
$\partial_{Y_{i}}$ act on $\polyalg[\para]$ 
as \emph{formal derivations}, {\ie} for $r\in\Z$ and $s\in\N$ we define
\begin{gather*}
\partial_{X_{i}}X_{j}^{\lambda+r}=\delta_{i,j}(\lambda+r)
\cdot X_{i}^{\lambda+r-1}
,
\partial_{X_{i}}Y_{j}^{s}=0
,\\
\partial_{Y_{i}}Y_{j}^{s}=\delta_{i,j}s\cdot Y_{i}^{s-1}
,
\partial_{Y_{i}}X_{j}^{\lambda+r}=0
,
\end{gather*}
and we then extend these 
rules to all of $\polyalg[\para]$ linearly and by the Leibniz rule.
\end{Definition}

We let the algebra $\ugl[2]$ 
act on $\polyalg[\para]$ by 
\begin{gather}\label{Eq:ProofDiffOp}
E\mapsto e=X\partial_{Y}
,\quad
F\mapsto f=Y\partial_{X}
,\quad
L_{1}\mapsto l_{1}=X\partial_{X}
,\quad
L_{2}\mapsto l_{2}=Y\partial_{Y}.
\end{gather}
The action \autoref{Eq:ProofDiffOp} extends to an action of 
$\ugl[2]$ on 
\begin{gather*}
\polyalg\cong\medotimes_{i=1}^{n}\polyalg[\para_{i}]
\end{gather*}
by using the 
usual coproduct of $\ugl[2]$ determined by $\Delta(x)=x\otimes 1+1\otimes x$ 
for all $x\in\gl[2]$.

We have a dual action of $\ugl$ on $\polyalg$ determined by 
\begin{gather}\label{Eq:ProofDiffOpTwo}
E_{ij}\mapsto e_{ij}=X_{i}\partial_{X_{j}}+Y_{i}\partial_{Y_{j}}.
\end{gather}
In particular, $E_{i}$ acts as $e_{i(i+1)}$, 
$F_{i}$ acts as $e_{(i+1)i}$ 
and $L_{i}$ acts as $e_{ii}$.

For $r\in\Z$
and $s\in\N$ let $X_{i}^{\lambda+r}Y_{i}^{s}$ be of degree $r+s$.
This gives us a $\Z^{n}$ grading on $\polyalg$.
For $\bsym[d]=(d_{1},\dots,d_{n})\in\Z^{n}$ 
we denote the $\Z^{n}$ graded piece
of $\polyalg$ of degree 
$\bsym[d]$ by $(\polyalg)_{\bsym[d]}$.

\begin{Lemma}\label{L:ProofDiffOp}
The graded $\K$ vector space
\begin{gather*}
\polyalg\cong\medoplus_{\bsym[d]\in\Z^{n}}(\polyalg)_{\bsym[d]}
\end{gather*} 
is an $\ugl[2]$ module when endowed with \autoref{Eq:ProofDiffOp} 
that is isomorphic to $\verma{\oplus\bsym}$ that decomposes as above.
Moreover, it is also an $\ugl$ module when endowed with \autoref{Eq:ProofDiffOpTwo},
and the two actions commute.
\end{Lemma}

\begin{proof}
That \autoref{Eq:ProofDiffOp} defines a 
homogeneous action of $\ugl[2]$ 
is easy to see. 

The resulting $\ugl[2]$ module 
is isomorphic to $\verma{\oplus\bsym}$ as in the 
classical $\gl[2]$ theory: For $d=0$ the basis elements 
of $(\polyalg[\para])_{0}$ are of the form $X^{\para-r}Y^{r}$ 
for $r\in\N$ and {\eg} $f(X^{\para-r}Y^{r})
=(\para-r)\cdot X^{\para-r-1}Y^{r+1}$. Comparing this 
with the classical version of \autoref{E:DualityVerma} shows that 
$(\polyalg[\para])_{0}\cong\verma{\para}$. For general $d\in\Z$ 
the story is just shifted and we get 
$(\polyalg[\para])_{d}\cong\verma{\para+d}$. 
These isomorphisms
extend to $\polyalg\cong\verma{\oplus\bsym}$ by using the coproduct.

That \autoref{Eq:ProofDiffOpTwo} defines 
an $\ugl[n]$ action and that the two actions commute
are direct calculations.
\end{proof}

We always use the two actions \autoref{Eq:ProofDiffOp} 
and \autoref{Eq:ProofDiffOpTwo} for the remainder of this section.
Note that
\autoref{L:ProofDiffOp} gives us a 
$\ugl[2]$-$\ugl[n]$ bimodule structure on $\polyalg$.

\begin{Notation}\label{N:ProofPowerNotation}
For $\bsym[Z]=(Z_{1},\dots,Z_{n})$ and 
$\bsym[b]=(b_{1},\dots,b_{n})$ we write 
$\bsym[Z]^{\bsym[b]}=Z_{1}^{b_{1}}\cdot\dots\cdot
Z_{n}^{b_{n}}$.
\end{Notation}

\begin{Lemma}\label{L:ProofHW}
The element $X^{\para+r}Y^{s}$ is 
annihilated by $E\in\gl[2]$ if and only if $s=0$.
\end{Lemma}

\begin{proof}
This holds since $e=X\partial_{Y}$ so that 
$e(X^{\para+r}Y^{s})=s\cdot X^{\para+r+1}Y^{s-1}$.
\end{proof}

An $\ugl[k]$ module is called \emph{countable semisimple} 
if it is a countable direct sum of countable dimensional 
simple $\ugl[k]$ modules.

\begin{Lemma}\label{L:ProofSemisimple}
For admissible parameters $\bsym$
the $\ugl[2]$ module $\polyalg$ is countable semisimple.
\end{Lemma}

\begin{proof}
Extending {\eg} 
\cite[Section 2]{Ka-tensor-infinite-dimensional-modules} to $\Kpara$, 
we let $\tilde{\mathcal{O}}$ denote enlarged 
category $\mathcal{O}$. We will not define $\tilde{\mathcal{O}}$ 
here as it can be defined, {\muta}, as in \cite[Section 2]{Ka-tensor-infinite-dimensional-modules}
with the same properties as therein. In particular, we have $\polyalg\cong\verma{\oplus\bsym}\in\tilde{\mathcal{O}}$. 

By the usual Yoga,  
the $\verma{\para_{i}}$ are costandard objects in 
$\tilde{\mathcal{O}}$. The usual Yoga, see e.g. \cite[Proposition 2.7]{Ka-tensor-infinite-dimensional-modules}
or \cite[Section 2]{AnStTu-cellular-tilting} 
and the extra notes for that paper in the arXiv version of it,
also gives that tensor products of 
$\ugl[2]$ modules with a costandard filtration 
have a costandard filtration. Moreover, the condition 
$\sums\para_{i}\notin\Z$ for all $i\in\{1,\dots,n\}$ ensures that 
all appearing costandard filtration factors have highest weight not being 
in $\Z$. Thus, all costandard filtration factors are tilting 
since they are simple and costandard, which is also a consequence of the usual Yoga.

It then follows
that $\verma{\oplus\bsym}$ decomposes into a direct sum 
of indecomposable tilting $\ugl[2]$ modules in $\tilde{\mathcal{O}}$,
and tracking the highest weight as in 
\autoref{L:ProofHW} and using that $\bsym$ 
is admissible shows that these 
indecomposable tilting $\ugl[2]$ modules are actually simple and
of the form $\verma{\mu}$ for generic $\mu\in\K$.

Finally, everything involved is clearly countable, so we are done.
\end{proof}

\begin{Lemma}\label{L:ProofHWTwo}
Let $\bsym$ be admissible.
As $\ugl[2]$ modules we have
\begin{gather}\label{Eq:ProofDecomTensor}
\polyalg\cong\verma{\oplus\bsym}\cong
\medoplus_{\substack{\countingpara\in\Z\\ t\in\N}}
\verma{\sums{\para_{n}}+\countingpara-t,t}\otimes\densemult{\sums{\para_{n}}+\countingpara-t,t}
,
\end{gather}
where $\densemult{\sums{\para_{n}}+\countingpara-t,t}$ 
is a multiplicity $\K$ vector space.
\end{Lemma}

\begin{proof}
For generic $\para\in\K$ and any $\para^{\prime}\in\K$
one can decompose $\verma{\para}\otimes\verma{\para^{\prime}}$ 
explicitly, {\ie}
as $\ugl[2]$ modules we have
\begin{gather}\label{Eq:ProofDecomTensorTwo}
\verma{\para}\otimes\verma{\para^{\prime}}
\cong
\medoplus_{t\in\N}
\verma{\para+\para^{\prime}-t,t}\otimes
\densemult{\para+\para^{\prime}-t,t},
\end{gather}
for some countable dimensional 
multiplicity $\K$ vector space $\densemult{\para+\para^{\prime}-t,t}$.
This decomposition \autoref{Eq:ProofDecomTensorTwo} 
follows from \autoref{L:ProofHW}
and \autoref{L:ProofSemisimple}
and the universal property of Verma modules.

More general, the decomposition \autoref{Eq:ProofDecomTensorTwo} 
can then be proven by using the proof of \autoref{L:ProofSemisimple}
which shows that $\polyalg$ is a(n infinite) direct sum of 
of simple tilting $\ugl[2]$ modules. The point is that the characters 
of simple tilting $\ugl[2]$ modules are well-known, since these are Verma 
modules, and we 
of course know the character of $\polyalg$ itself. Using this and semisimplicity \autoref{L:ProofSemisimple}, we hence get the 
claimed formula by successively identifying the characters in 
$\polyalg$. That is, we first get
\begin{gather*}
\polyalg\cong\verma{\oplus\bsym}\cong
\medoplus_{\bsym[d]\in\Z^{n}}
\verma{\para_{1}+d_{1}}
\otimes\dots\otimes\verma{\para_{n}+d_{n}}
\cong
\medoplus_{\substack{\bsym[d]\in\Z^{n}\\ t\in\N}}
\verma{\sums{\para_{n}}+\sums{d_{n}}-t,t}\otimes\densemult{\sums{\para_{n}}+\sums{d_{n}}-t,t}
,
\end{gather*}
and then grouping isomorphic 
$\ugl[2]$ modules gives \autoref{Eq:ProofDecomTensor}.
\end{proof}

We now aim to identify 
$\densemult{\sums{\para_{n}}+\countingpara-t,t}$ 
from \autoref{Eq:ProofDecomTensor} explicitly.
To this end, we define a $\K$ sub vector space $\moduledet[{b,c}]$ 
of $\polyalg$ that we will use for this purpose:

\begin{Definition}\label{D:ProofBasisOne}
Write $\detsym=\det
\begin{psmallmatrix}x_{i}&y_{i}\\x_{j}&y_{j}\end{psmallmatrix}$,
$\detsym[i]=\detsym[{i(i+1)}]$, $\bsym[a]=(\detsym[{1}],\dots,\detsym[{n-1}])$ 
and $\bsym[l]=(l_{1}\dots,l_{n-1})\in\N^{n-1}$.	
For $b\in\Z$ and $c\in\N$ let
\begin{gather*}
\moduledet[{b,c}]=\K\basisdet\subset\polyalg
,
\end{gather*}
where $\basisdet=\basisdet(b,c)=\{\bsym[X]^{\bsym[\para]+\bsym[r]}
\bsym[a]^{\bsym[l]}|\sums{r_{n}}=b-c,\sums{l_{n-1}}=c\}$.
\end{Definition}

\begin{Lemma}\label{L:ProofBasesOne}
For fixed $c\in\N$	
let $\medprod_{r,s}\detsym[{rs}]$ be 
a product of $c$ determinants. Then $\medprod_{r,s}\detsym[{rs}]
\in\K[\bsym[X]^{\pm 1}]\{\bsym[a]^{\bsym[l]}|\sums{l_{n-1}}=c\}$.
\end{Lemma}

\begin{proof}
The determinant of the singular matrix
\begin{gather*}
\det\begin{psmallmatrix}X_{r} & X_{r} & Y_{r}\\X_{i} & X_{i} & Y_{i}\\X_{s} & X_{s} & Y_{s}\end{psmallmatrix}
=
X_{r}\detsym[{is}]-X_{i}\detsym[{rs}]+X_{s}\detsym[{ri}]=0
\end{gather*}
gives the relation $\detsym[{rs}]=X_{i}^{-1}(X_{r}\detsym[{is}]+X_{s}\detsym[{ri}])$ 
for all $i\in\{1,\dots,n\}$. This relation can 
then be successively applied to prove the statement.
\end{proof}

\begin{Lemma}\label{L:ProofBases}
The $\ugl$ action from \autoref{L:ProofDiffOp}
stabilizes $\moduledet[{b,c}]\subset\polyalg$.
\end{Lemma}

\begin{proof}
A straightforward calculation gives
\begin{gather*}
e_{ij}(\detsym[{rs}])=
\begin{cases}
\detsym[{is}]&\text{if }j=r,
\\
\detsym[{ri}]&\text{if }j=s,
\\
0&\text{else}.
\end{cases}
\end{gather*}
Using this we get
\begin{gather*}
e_{ij}(\bsym[X]^{\bsym[\para]+\bsym[r]}
\bsym[a]^{\bsym[l]})
=
e_{ij}(\bsym[X]^{\bsym[\para]+\bsym[r]})\bsym[a]^{\bsym[l]}
+
\bsym[X]^{\bsym[\para]+\bsym[r]}e_{ij}(\bsym[a]^{\bsym[l]})
\\
\scalebox{0.9}{$=(\para_{i}+r_{i})\cdot
\bsym[X]^{\bsym[\para]+\bsym[r]+\stdvec[i]-\stdvec[j]}
\bsym[a]^{\bsym[l]}
+l_{j-1}\cdot
\bsym[X]^{\bsym[\para]+\bsym[r]}\bsym[a]^{\bsym[l]-\stdvec[j-1]}
\underbrace{\detsym[{(j-1)i}]}_{\text{rewrite}}
+l_{j}\cdot
\bsym[X]^{\bsym[\para]+\bsym[r]}\bsym[a]^{\bsym[l]-\stdvec[j]}
\underbrace{\detsym[{i(j+1)}]}_{\text{rewrite}}.$}
\end{gather*}
Now we use the 
rewriting as in the proof of \autoref{L:ProofBasesOne} on 
the marked terms, and we are done. Explicitly, we get
\begin{gather}\label{Eq:ProofFormulas}
\begin{aligned}
e_{ii}(\bsym[X]^{\bsym[\para]+\bsym[r]}
\bsym[a]^{\bsym[l]})
=&(\para_{i}+r_{i}+l_{i-1}+l_{i})\cdot
\bsym[X]^{\bsym[\para]+\bsym[r]}
\bsym[a]^{\bsym[l]}
,
\\
e_{i(i+1)}(\bsym[X]^{\bsym[\para]+\bsym[r]}
\bsym[a]^{\bsym[l]})
=&
(\para_{i}+r_{i}+l_{i+1})\cdot
\bsym[X]^{\bsym[\para]+\bsym[r]+\sroot[i]}
\bsym[a]^{\bsym[l]}
+l_{i+1}\cdot
\bsym[X]^{\bsym[\para]+\bsym[r]-\sroot[i+1]}\bsym[a]^{\bsym[l]+\sroot[i]}
,
\\
e_{(i+1)i}(\bsym[X]^{\bsym[\para]+\bsym[r]}
\bsym[a]^{\bsym[l]})
=&
(\para_{i}+r_{i}+l_{i-1})\cdot
\bsym[X]^{\bsym[\para]+\bsym[r]-\sroot[i]}
\bsym[a]^{\bsym[l]}
+l_{i-1}\cdot
\bsym[X]^{\bsym[\para]+\bsym[r]+\sroot[i-1]}\bsym[a]^{\bsym[l]-\sroot[i]}.
\end{aligned}
\end{gather}
Here we used
$\detsym[{i(i+2)}]=X_{i+1}^{-1}(X_{i}
\detsym[{i+1}]+X_{i+2}\detsym[{i}])$
as well as $\detsym[{(i-1)(i+1)}]=X_{i}^{-1}(X_{i-1}
\detsym[{i}]+X_{i+1}\detsym[{i-1}])$.
\end{proof}

We want to show that $\moduledet[{b,c}]$ is a dense module 
as in \autoref{T:DualityHowe}.
To do this we need an analog of the GT basis,  
and to define it we need to prepare
the definition with some preliminaries.

\begin{Notation}\label{N:ProofGTBasis}
\leavevmode

\begin{enumerate}

\item For $\bsym[d]=(d_{1},\dots,d_{n-1})\in\N^{n-1}$, 
we denote by $\binom{\bsym[d]}{\bsym[s]}\in\N$ the 
\emph{multino\-mial-type} number
defined by the expansion
$\medprod_{i=1}^{n-1}(\sums{X_{i}})^{d_{i}}=
\medsum_{\bsym[s]}\binom{\bsym[d]}{\bsym[s]}\cdot\bsym[X]^{\bsym[s]}$. 

\item We
let $(k)_{l}=k(k+1)\cdot\dots\cdot(k+l-1)$ 
be the (increasing) \emph{Pochhammer symbol}.

\item We write $\prods{(\bsym[\para],\bsym[r],\bsym[d],\bsym[j])}$ 
for
$(\sums{\para_{1}}+\sums{r_{1}}-j_{1}+1)_{d_{1}+j_{1}-j_{0}}\cdot
\dots\cdot(\sums{\para_{n-1}}+\sums{r_{n-1}}-j_{n-1}+1)_{d_{n-1}+j_{n-1}-j_{n-2}}$, 
where $\bsym[j]\in\Z^{n}$ with $j_{n-1}=j_{n}=0$.

\end{enumerate}
\end{Notation}

\begin{Definition}\label{D:ProofGTBasis}
For the GT pattern
\begin{gather*}
\gt[{\bsym[x],\bsym[c]}]=
\begin{array}{CCCCCCCCC}
\cellcolor{tomato!50}$x_{n}$ & & \cellcolor{spinach!50}$c_{n}$ & \cellcolor{spinach!50} & \cellcolor{spinach!50}$0$ & \cellcolor{spinach!50} & \cellcolor{spinach!50}$\dots$ & \cellcolor{spinach!50} & \cellcolor{spinach!50}$0$ \\
& \cellcolor{tomato!50}$x_{n{-}1}$ & & \cellcolor{spinach!50}$c_{n{-}1}$ & \cellcolor{spinach!50} & \cellcolor{spinach!50}$\ddots$ & \cellcolor{spinach!50} & \cellcolor{spinach!50}$\rddots$ & \\
& & \cellcolor{tomato!50}$\ddots$ & & \cellcolor{spinach!50}$\ddots$ & \cellcolor{spinach!50} & \cellcolor{spinach!50}$0$ & & \\
& & & \cellcolor{tomato!50}$x_{2}$ & & \cellcolor{spinach!50}$c_{2}$ & & & \\
& & & & \cellcolor{tomato!50}$x_{1}$ & & & &
\end{array}
,\\
\begin{gathered}
c_{1}=0,c_{n}=c,
\\
\bsym[d]=(c_{2}-c_{1},c_{3}-c_{2},\dots,c_{n}-c_{n-1}),
\\
\sums{r_{i}}=
x_{i}-\sums{\para_{i}}-\sums{c_{i-1}}+c_{i+1}
,
\end{gathered}
\end{gather*}
we define $\bsym[d]=(d_{1},\dots,d_{n-1})$ and 
$\bsym[r]=(r_{1},\dots,r_{n})$ as above. The associated 
\emph{GT vector} is
\begin{gather}\label{Eq:ProofGTVectors}
\scalebox{0.85}{$\gt[{\bsym[d],\bsym[r]}]
=\medsum_{\bsym[j]\in\N^{n-2}}
\medbinom{\bsym[d]}{\bsym[d]+\sum_{i=1}^{n-2}j_{i}\sroot[i]}
\prods{(\bsym[\para],\bsym[r],\bsym[d],\bsym[j])}
\cdot
\bsym[X]^{\bsym[\para]+\bsym[r]-\sum_{i=1}^{n}(j_{i}-j_{i-2})\stdvec[i]}
\bsym[a]^{\bsym[d]+\sum_{i=1}^{n-2}j_{i}\sroot[i]}\in\moduledet[{b,c}].$}
\end{gather}
The set of these GT vectors is denoted by $\basisgt$.
\end{Definition}

The following manipulation of one of the scalars 
defining $\gt[{\bsym[d],\bsym[r]}]$ will come in handy.

\begin{Lemma}\label{L:ProofGTBasis}
We have $\medbinom{\bsym[d]}{\bsym[d]+\sum_{i=1}^{n-2}j_{i}\sroot[i]}=\prod_{i=1}^{n-2}\binom{d_{i+1}+j_{i+1}}{j_{i}}$.
\end{Lemma}

\begin{proof}
By using
$(\sums{X_{i}})^{d_{i}}=(\sums{X_{i-1}}+X_{i})^{d_{i}}=
\sum_{j=0}^{d_{i}}\binom{d_{i}}{j}(\sums{X_{i-1}})^{d_{i}-j}X_{i}^{j}$	
recursively.
\end{proof}

Let $z_{k}=\medsum_{1\leq j<i\leq k}e_{ji}e_{ij}$ 
for $k\in\{1,\dots,n\}$.
We need what we call the \emph{Casimir elements} of $\gl[n]$, which are 
defined for $k\in\{1,\dots,n\}$ as
\begin{gather}\label{Eq:ProofCasimirElement}
\casimir[k]
=\medsum_{1\leq i \neq j \leq k}E_{ij}E_{ji}+\medsum_{i=1}^{k}E_{ii}^{2}
=2z_{k}+\medsum_{1\leq j<i\leq k}(E_{ii}-E_{jj})+\medsum_{i=1}^{k}E_{ii}^{2}
.
\end{gather}
The notation is such that $\casimir[k]$ is the 
usual Casimir element of $\gl[k]$. We write $\casimirsmall[k]$ 
for the associated operator.
For \autoref{L:ProofCasimir} below, which is 
the main lemma regarding the Casimir 
elements, we need the following formula for the action of
$z_{k}$.

\begin{Lemma}\label{L:ProofFormulas}
We have $z_{k}(\bsym[X]^{\bsym[\para]+\bsym[r]}
\bsym[a]^{\bsym[l]})=s\cdot\bsym[X]^{\bsym[\para]+\bsym[r]}
\bsym[a]^{\bsym[l]}+e$ with $s\in\K$ and an error term 
$e$ given by
\begin{align*}
s=
&\medsum_{i=1}^{k}
\Big(
(\para_{i}+r_{i})(\sums{\para_{i-1}}+\sums{r_{i-1}}
+\sums{l_{n}}-l_{i-1}-l_{i}+i-1)
\\
&+(i-1)l_{i}+(i-2)l_{i-2}+(l_{i-1}+l_{i})\sums{l_{i-2}}\Big)
,
\\
e=&
l_{k}\medsum_{i=1}^{k-1}(\para_{i}+r_{i})\cdot
\bsym[X]^{\bsym[\para]+\bsym[r]+\stdvec[k]+\stdvec[{k+1}]-\stdvec[i]-\stdvec[{i+1}]}\bsym[a]^{\bsym[l]-\stdvec[k]+\stdvec[i]}
.
\end{align*}
\end{Lemma}

\begin{proof}
A tedious calculation using the previous formulas.
\end{proof}

\begin{Lemma}\label{L:ProofCasimir}
Let $\bsym$ be admissible. 
The Casimir elements separate $\basisgt$ (on weight spaces), and $\basisgt$ is a basis of $\moduledet[{b,c}]$.
\end{Lemma}

\begin{proof}
The proof splits into three steps.

\textit{Separation.}
We first assume that the Casimir elements act by a scalar on 
$\basisgt$. We let 
$\nu_{k}=x_{k}\stdvec[1]+c_{k}\stdvec[2]$ where we recall that 
$x_{k}=\sums{r_{k}}-\sums{\para_{k}}-\sums{c_{k-1}}+c_{k+1}$.
We assume the scalar is
\begin{gather}\label{Eq:ProofCasimir}
\casimir[k]\text{ acts on $\gt[{\bsym[d],\bsym[r]}]$ by }
\langle
\nu_{k}+2\rho^{(k)},\nu_{k}
\rangle
=x_{k}(x_{k}+k-1)+c_{k}(c_{k}+k-3)
,
\end{gather}
where $\rho^{(k)}=\frac{1}{2}\medsum_{i=1}^{k}(k-2i+1)\cdot\stdvec[i]$ 
mimics the usual half-sum of the positive roots of $\gl[n]$.

On a weight space we have $c_{k}=a-x_{k}$ for some $a\in\Z$. 
Hence, we get the parabola $2x_{k}^{2}+(2-2a)x_{k}+a^{2}+ka-3a$ 
from \autoref{Eq:ProofCasimir}. Assume that there are two values 
$x_{k}$ and $x_{k}^{\prime}$ as in the nonintegral part of GT patterns which satisfy this parabola. Solving $2x_{k}^{2}+(2-2a)x_{k}+a^{2}+ka-3a=2(x_{k}^{\prime})^{2}+(2-2a)x_{k}^{\prime}+a^{2}+ka-3a$ gives either $x_{k}=x_{k}^{\prime}$ or $x_{k}+x_{k}^{\prime}-a=-1$. The second solution gives $x_{k}^{\prime}\in\Z$, which contradicts admissibility.

Note that this implies that the Casimir elements separate, 
so it remains to verify \autoref{Eq:ProofCasimir}.

\textit{Scalar verification.} We thus need to compute $\casimirsmall[k](\gt[{\bsym[d],\bsym[r]}])$.
The calculation that 
$\casimirsmall[k](\gt[{\bsym[d],\bsym[r]}])$ equals \autoref{Eq:ProofCasimir}
boils down to a longish manipulation of symbols where one reindexes the sum 
defining $\gt[{\bsym[d],\bsym[r]}]$ appropriately. We sketch the main 
step in this calculation now.
\begin{enumerate}[label=(\alph*)]

\item First, 
we use the second expression of $\casimir[k]$ in \autoref{Eq:ProofCasimirElement}. As before, we use  
$z_{k}=\medsum_{1\leq j<i\leq k}e_{ji}e_{ij}$ and 
we also write $h_{k}$ for the Cartan part so that 
$\casimirsmall[k]=2z_{k}+h_{k}$.
By \autoref{Eq:ProofFormulas},
the Cartan part gives 
$h_{k}(\bsym[X]^{\bsym[\para]+\bsym[r]}
\bsym[a]^{\bsym[l]})=s^{\prime}\cdot\bsym[X]^{\bsym[\para]+\bsym[r]}
\bsym[a]^{\bsym[l]}$ with scalar
\begin{gather*}
s^{\prime}
=\medsum_{1\leq j<i\leq k}
(\para_{i}+r_{i}+l_{i-1}+l_{i}-\para_{j}-r_{j}-l_{j-1}-l_{j})
\\+\medsum_{i=1}^{k}(\para_{i}+r_{i}+l_{i-1}+l_{i})^{2}\in\K.
\end{gather*}

\item We also have the scalar $s$ from \autoref{L:ProofFormulas}. 
Thus, we get, again using \autoref{L:ProofFormulas}, that
\begin{gather}\label{Eq:ProofAnotherFormula}
\casimirsmall[k](\bsym[X]^{\bsym[\para]+\bsym[r]}
\bsym[a]^{\bsym[l]})
=
(2s+s^{\prime})
\cdot
\bsym[X]^{\bsym[\para]+\bsym[r]}
\bsym[a]^{\bsym[l]}
+e,
\end{gather}
where $e$ is the error term in \autoref{L:ProofFormulas}.

\item Next, we need to take the sum 
of \autoref{Eq:ProofAnotherFormula} as in the definition of the GT vectors.
The resulting expression can then by manipulated as in the next few bullet points.	

\item We change the summation and use a few tricks 
to get the same expressions defining the GT vectors from 
\autoref{Eq:ProofGTVectors}:
\begin{enumerate}[label=$\bullet$]

\item For the part with the multinomial scalars 
we use \autoref{L:ProofGTBasis} and the well-known formula
$\binom{a-1}{b-1}=\frac{b}{a}\binom{a}{b}$ to rewrite {\eg}
\begin{gather*}
\binom{d_{i-1}+j_{i-1}\fcolorbox{black}{orchid!50}{$-1$}}{j_{i-2}\fcolorbox{black}{orchid!50}{$-1$}}
=
\frac{j_{i-2}}{d_{i-1}+j_{i-1}}
\binom{d_{i-1}+j_{i-1}\fcolorbox{black}{orchid!50}{$+0$}}{j_{i-2}\fcolorbox{black}{orchid!50}{$+0$}}
.
\end{gather*}
We further use $\binom{a-1}{b}=\frac{a-b}{a}\binom{a}{b}$ to rewrite for example
\begin{gather*}
\binom{d_{i-1}+j_{i-1}\fcolorbox{black}{orchid!50}{$-1$}}{j_{i-2}}
=
\frac{d_{i-1}+j_{i-1}-j_{i-2}}{d_{i-1}+j_{i-1}}
\binom{d_{i-1}+j_{i-1}\fcolorbox{black}{orchid!50}{$+0$}}{j_{i-2}}
\end{gather*}
Here we marked the parts that we change to match the GT vectors.
(If $a=0$ in these formulas, then we would use 
$a\binom{a-1}{b-1}=b\binom{a}{b}$ and 
$a\binom{a-1}{b}=(a-b)\binom{a}{b}$ which give 
$0=0$ so we can ignore these cases.)

\item We rewrite the Pochhammer symbols as well, for example:
\begin{gather*}
(\sums{\para}_{i}+\sums{r}_{i}-j_{i}\fcolorbox{black}{orchid!50}{$+2$})_{d_{i}+j_{i}-j_{i-1}\fcolorbox{black}{orchid!50}{$-1$}}
\\=
\frac{1}{\sums{\para}_{i}+\sums{r}_{i}-j_{i}+1}
(\sums{\para}_{i}+\sums{r}_{i}-j_{i}\fcolorbox{black}{orchid!50}{$+1$})_{d_{i}+j_{i}-j_{i-1}\fcolorbox{black}{orchid!50}{$+0$}}
.
\end{gather*}
We again highlight the parts we change to match the expressions in the GT vectors.

\end{enumerate}

\item All introduced fractions disappear in the end. To elaborate, we get
\begin{gather*}
\medsum_{i=1}^{k-1}
\bigg(\medprod_{l=i}^{k-1}\frac{j_{l}}{d_{l}+j_{l}}\bigg)
(d_{i}+j_{i}-j_{i-1})\\
\cdot
\frac{\medprod_{l=i+1}^{k}(\sums{\para_{l}}+\sums{r_{l}}+d_{l}-j_{l-1}+1)}
{\medprod_{l=i}^{k-1}(\sums{\para_{l}}+\sums{r_{l}}-j_{l}+1)}
(\sums{\para_{i}}+\sums{r}_{i}-j_{i}-j_{i-1}+1).
\end{gather*}
The sum of the two first terms is:
\begin{gather*}
\frac{\medprod_{i=2}^{k-1}j_{i}}{\medprod_{i=3}^{k-1}(d_{i}+j_{i})}
\frac{\medprod_{i=3}^{k}(\sums{\para_{i}}+\sums{r_{i}}+d_{i}-j_{i-1}+1)}{\medprod_{i=3}^{k-1}(\sums{\para_{i}}+\sums{r_{i}}-j_{i}+1)}
.
\end{gather*}
We continue, analyzing three, four {\etc} terms, until we find
\begin{gather*}
j_{k-1}(\sums{\para_{k}}+\sums{r_{k}}+d_{k}-j_{k-1}+1).
\end{gather*}

\item This implies that the overall scalar for the 
$\bsym[j]\in\N^{n-2}$ summand of $\gt[{\bsym[d],\bsym[r]}]$ is
\begin{gather}\label{Eq:ProofJCancles}
\begin{gathered}
\medsum_{i=3}^{k}(\para_{i}+r_{i}+d_{i}+d_{i-1})
(\sums{d_{i-2}}\fcolorbox{black}{blue!50}{$+j_{i-2}$})
\\+
\medsum_{i=2}^{k}\big(
(\sums{\para_{i-1}}+\sums{r_{i-1}}\fcolorbox{black}{blue!50}{$-j_{i-1}-j_{i-2}$})
(d_{i}\fcolorbox{black}{blue!50}{$+j_{i}-j_{i-1}$})
\\
+
(\para_{i}+r_{i}+d_{i}\fcolorbox{black}{blue!50}{$+j_{i-2}-j_{i-1}$})
(\sums{\para_{i-1}}+\sums{r_{i-1}}\fcolorbox{black}{blue!50}{$-j_{i-1}-j_{i-2}$}+i-1)
\\
+(i-2)(d_{i-1}\fcolorbox{black}{blue!50}{$+j_{i-1}-j_{i-2}$})
\big)
\fcolorbox{black}{blue!50}{$+j_{k-1}(\sums{\para_{k}}+\sums{r_{k}}+d_{k}-j_{k-1}+1)$}
.
\end{gathered}
\end{gather}
We marked the dependencies on $\bsym[j]$.

\item The dependencies on $\bsym[j]$ in 
\autoref{Eq:ProofJCancles} cancel and
we get $\casimirsmall[k](\gt[{\bsym[d],\bsym[r]}])=
s^{\prime\prime}\cdot\gt[{\bsym[d],\bsym[r]}]$ for the scalar
\begin{gather*}
s^{\prime\prime}=
\medsum_{i=3}^{k}(\para_{i}+r_{i}+d_{i}+d_{i-1})\sums{d_{i-2}}
+
\medsum_{i=2}^{k}\big(
(\sums{\para_{i-1}}+\sums{r_{i-1}})d_{i}
\\
+(\para_{i}+r_{i}+d_{i})
(\sums{\para_{i-1}}+\sums{r_{i-1}}+i-1)
+(i-2)d_{i-1}\big)
\\
+\medsum_{i=1}^{k}\big((k-2i+1)(\para_{i}+r_{i}+d_{i-1}+d_{i})
+(\para_{i}+r_{i}+d_{i-1}+d_{i})^{2}\big)
\in\K.
\end{gather*}

\item Finally, matching $s^{\prime\prime}$ with \autoref{Eq:ProofCasimir} is 
done by comparing the linear and the quadratic terms separately. 
This is again tedious, but straightforward.

\end{enumerate}
\textit{Basis.}
The linear independence of $\basisgt$ follows from
\autoref{Eq:ProofCasimir} and that $\basisgt$ spans 
follows because the definition of $\gt[{\bsym[d],\bsym[r]}]$ 
implies that $\basisgt$ is upper triangular
(with an appropriate order) to $\basisdet$.
\end{proof}

\begin{Remark}\label{R:ProofCasimir}
We were able to guess
the formulas in \autoref{Eq:ProofCasimir} 
because of significant help of Magma and Mathematica, 
which were used to find the GT bases expressions 
in \autoref{D:ProofGTBasis}, as well as 
the formula given in \cite[(6.2)]{MiToLa-casimir-braids}. 
The computationally expensive proof of \autoref{L:ProofCasimir} 
was then also obtained by computer help.
We however stress that everything can be done by hand and computers 
were only used to guess the various steps.
\end{Remark}

\begin{Lemma}\label{L:ProofIsSimpleTwo}
For admissible parameters $\bsym$ we 
have that $\moduledet[{b,c}]$ is a simple dense $\ugl[n]$ module that has a GT pattern realization.
\end{Lemma}

\begin{proof}
We first show that $\moduledet[{b,c}]$ is a simple 
$\ugl[n]$ module. To this end, we use that the Casimir elements separate
the GT patterns in $\basisgt$, see 
\autoref{L:ProofCasimir}, and then we use similar arguments as 
in \cite[Lemma 3]{Ma-gt-dense}.	That is, 
we claim that the $e_{ij}$ act injectively (and thus, bijectively) 
on $\gt[{\bsym[d],\bsym[r]}]$
and also that the action graph 
of the $e_{ij}$ action on $\basisgt$ is strongly 
connected. 

The first claim follows from 
\autoref{L:ProofCasimir}, which implies 
that it is enough to show injectivity on $\basisgt$, 
and the formulas 
for the action of
$e_{ij}$ on $\gt[{\bsym[d],\bsym[r]}]$ 
that we get from \autoref{Eq:ProofFormulas}.

For the second claim we compute that
\begin{align*}
e_{i(i+1)}(\gt[{\bsym[d],\bsym[r]}])=&
\frac{\left(\sums{\para_{i}}+\sums{r_{i}}+1\right)
\left(\sums{\para_{i+1}}+\sums{r_{i+1}}+d_{i+1}\right)}
{\sums{\para_{i}}+\sums{r_{i}}+d_{i}+1}
\cdot
\gt[{\bsym[d],\bsym[r]+\sroot[i]}]
\\&+
\frac{d_{i}\left(\sums{\para_{i+1}}+\sums{r_{i+1}}+d_{i}+d_{i+1}+1\right)}
{\sums{\para_{i}}+\sums{r_{i}}+d_{i}+1}
\cdot\gt[{\bsym[d]+\sroot[i-1],\bsym[r]-\sroot[i-1]}],
\end{align*}
with the second term being zero if $i=1$ or if $\bsym[d]+\sroot[i-1]\not\in\N^{n-1}$. There is also
a similar formula for $e_{(i+1)i}(\gt[{\bsym[d],\bsym[r]}])$ 
with swapped signs in front of the $\sroot[j]$ and similar coefficients.
Note that all appearing coefficients are nonzero since we 
have admissible parameters.

Thus, the action graph is strongly connected and 
hence, the $e_{ij}$ act bijectively and have a strongly 
connected action graph, showing that $\moduledet[{b,c}]$ is a simple 
$\gl[n]$ module.

Finally, it follows from the definitions that $\moduledet[{b,c}]$ 
is a dense $\ugl$ module in the sense of {\eg} the introduction 
of \cite{Ma-gt-dense}.
\end{proof}

Let $\dense{b,c}$ denote the dense $\ugl[n]$ module 
defined in \cite[Section 3]{Ma-gt-dense} associated 
to a two diagonal GT pattern.

\begin{Proposition}\label{P:ProofMatch}
Assume that we have admissible parameters 
satisfying $\para_{i}\notin\Z$ for $i\in\{1,\dots,n\}$. We have an isomorphism of 
$\ugl$ modules $\moduledet[{b,c}]\cong\dense{b,c}$.
\end{Proposition}

\begin{proof}
We have also verified that 
$\moduledet[{b,c}]$ is simple and dense in 
\autoref{L:ProofIsSimpleTwo}. Thus, we can use the 
classification of these modules from \cite{Ma-classification-weight-modules}, 
see also \cite[Section 2.3]{Ma-gt-dense}.
\end{proof}

\begin{Remark}\label{R:ProofMatch}
Note that the isomorphism in \autoref{P:ProofMatch} is not explicit.
Any explicit isomorphism would divide, or multiply, by our parameters plus integers. 
That is why we need the assumption 
$\para_{i}\notin\Z$ for $i\in\{1,\dots,n\}$ in \autoref{P:ProofMatch}.
\end{Remark}

Abusing notation, we will write $\dense{b,c}$ instead of $\moduledet[{b,c}]$ to refer to its underlying GT pattern realization.

\begin{proof}[Proof of the classical version of \autoref{T:DualityHowe}]
There are three sta\-tements to verify.

\textit{Commuting actions.}
By \autoref{L:ProofDiffOp},
we can consider the $\ugl[2]$-$\ugl$ bimodule $\polyalg$ 
and it remains to verify the centralizer 
property and the $\ugl[2]$-$\ugl$ bimodule decomposition of \autoref{T:DualityHowe}.

All parameters in this proof are admissible from now on.

\textit{Bimodule decomposition.}
Let $b=\countingpara-t$ and $c=t$.
The $\ugl[2]$-$\ugl$ bimodule decomposition follows 
from \autoref{L:ProofIsSimpleTwo} after identifying 
$\moduledet[{b,c}]$ with the multiplicity space 
$\densemult{\sums{\para_{n}}+\countingpara-t,t}$ from \autoref{L:ProofHWTwo} 
as a $\K$ vector space.
Note that $\moduledet[{b,c}]$ is a 
$\ugl[n]$ module so it has a $\Z^{n}$ grading coming 
from the $\ugl[n]$ weight spaces. At the same time, because 
$\polyalg$ is an $\ugl[2]$-$\ugl$ bimodule, \autoref{L:ProofHWTwo} 
implies that $\densemult{\sums{\para_{n}}+\countingpara-t,t}$ is also a 
$\ugl[n]$ module, so we also have the notion of $\ugl[n]$ weight spaces.
Explicitly, $m$ in either $\moduledet[{b,c}]$ or 
$\densemult{\sums{\para_{n}}+\countingpara-t,t}$ is of degree $\bsym[d]\in\Z^{n}$ 
if $e_{ii}(m)=(\para_{i}+d_{i})\cdot m$.
We apply this definition and \autoref{Eq:ProofFormulas} 
to $\basisdet$ and get
\begin{gather*}
\grdim\moduledet[{b,c}]
=\medsum_{\bsym[d]\in\Z^{n}}
\delta_{b+c,\sums{d_{n}}}
\medbinom{c+n-2}{c}\bsym[Z]^{\bsym[d]},
\end{gather*}
where we use $\bsym[Z]=(Z_{1},\dots,Z_{n})$ to keep track of the 
graded pieces and $\grdim$ means graded dimensions. 
Moreover, using $\medsum_{c=0}^{t}
\binom{c+n-2}{c}=\binom{t+n-1}{t}$, we get
\begin{gather*}
\grdim\polyalg
=\medsum_{\substack{\bsym[d]\in\Z^{n}\\ t\in\N}}
\medbinom{t+n-1}{t}\bsym[Z]^{\bsym[d]}
=
\medsum_{\substack{\bsym[d]\in\Z^{n}\\ t\in\N}}
\bigg(
\medsum_{c=0}^{t}
\medbinom{c+n-2}{c}
\bigg)
\bsym[Z]^{\bsym[d]}
.
\end{gather*}
Thus, \autoref{Eq:ProofDecomTensor} implies that
\begin{gather}\label{Eq:ProofDenseGrDim}
\grdim\moduledet[{b,c}]
=\grdim\densemult{\sums{\para_{n}}+\countingpara-t,t}
=\medsum_{\bsym[d]\in\Z^{n}}
\delta_{b+c,\sums{d_{n}}}
\medbinom{c+n-2}{c}\bsym[Z]^{\bsym[d]}.
\end{gather}
Finally, note that $\moduledet[{b,c}]
\subset\densemult{\sums{\para_{n}}+\countingpara-t,t}$ because 
$e(\basisdet)=\{0\}$, as a simple calculation shows. Hence, 
$\moduledet[{b,c}]=\densemult{\sums{\para_{n}}+\countingpara-t,t}$.

\textit{Dense.}
Thus, it remains to prove that have 
dense subrings induced from the 
$\ugl[2]$ and $\ugl$ actions.
We denote their images in $\End_{\K}(\polyalg)$ 
by $\ring[S]$ and $\ring[T]$, respectively.
The $\ugl[2]$-$\ugl$ bimodule decomposition implies 
that $\ring[T]$ is dense in $\ring[S]^{\prime}$ 
since 
\begin{gather*}
\ugl[2]
\dgen
\End_{\ugl[n]}\big(
\medoplus_{\substack{\countingpara\in\Z\\ t\in\N}}
\verma{\sums{\para_{n}}+\countingpara-t,t}\otimes\dense{\countingpara-t,t}
\big)
\cong
\End_{\ugl[n]}(\polyalg)
.
\end{gather*}
Similarly, with swapped roles of 
$\ugl[2]$ and $\ugl$, 
we get that $\ring[S]$ is dense in $\ring[T]^{\prime}$.
\end{proof}

\begin{Remark}\label{R:ProofDenseWSpace}
Note that \autoref{Eq:ProofDenseGrDim} implies that 
the dense modules we use have weight spaces of 
constant dimension. This is actually true in 
more generality, see {\eg} \cite[Lemma 2]{Ma-gt-dense}.
\end{Remark}

%%%%%%%%%%%%%%%%%%%%%%%%%%%%%%%%%%%%%%%%%

\subsection{The quantum case}\label{SS:ProofQuantum}

%%%%%%%%%%%%%%%%%%%%%%%%%%%%%%%%%%%%%%%%%

We now specify our quantum conventions.

\begin{Notation}\label{N:ProofSpecialization}
\leavevmode
\begin{enumerate}

\item Let $\bsymv=(\parav_{1},\dots,\parav_{n})$ denote a tuple of variables.
Let $\aform=\Z[\vpar,\vpar^{-1}]$ 
denote the \emph{A-form}, $\afrac=\Q(\vpar)$ 
the field of fractions of $\aform$, $\apara=\aform[\bsymv,\vpar^{\parav_{1}},\dots,\vpar^{\parav_{n}}]$ 
and $\afracpara$ the field of fractions of $\apara$.

\item For a fixed ring $\ring[S]$ and choices $\bar{\vpar}$, 
$\bar{\vpar}^{-1}\in\ring[S]$ and $\bar{\bsymv}=
(\bar{\parav}_{1},\dots,\bar{\parav}_{n})\in\ring[S]^{n}$ 
such that $\bar{\vpar}^{\bar{\parav}_{1}},\dots,\bar{\vpar}^{\bar{\parav}_{n}}
\in\ring[S]$,
we can specialize any 
construction defined with coefficients from $\apara$ 
using $\placeholder\otimes_{\apara}\ring[S]$ where we see 
$\ring[S]$ as an $\apara$ module by $\vpar^{\pm 1}
\mapsto\bar{\vpar}^{\pm 1}$ and $\bsymv\mapsto\bar{\bsymv}$. 
We will use this always without the bar notation.
We can similarly specialize from $\aform$ instead of 
$\apara$ in the very same way. 

\item The \emph{classical} and the \emph{quantum specializations} are 
$\vpar\mapsto 1$ and $\vpar\mapsto\qpar$, respectively, and 
$\bsymv\mapsto\bsym$ for $\ring[S]=\K$ as fixed in \autoref{SS:DualityVermaDense}. Similarly for $\aform$ instead of 
$\apara$.

\end{enumerate}
\end{Notation}

For $k\in\Z_{\geq 1}$, let $\uagl[k]$ denote the 
\emph{quantum enveloping algebra} over $\aform$ of $\gl[k]$.
We use the conventions, excluding the Hopf algebra structure, from \cite{Lu-qgroups-root-of-1} 
or \cite{AnPoWe-representation-qalgebras} with $K_{i}^{\pm 1}=L_{i}^{\pm 1}L_{i+1}^{\mp 1}$. The $\aform$ algebra $\uagl[k]$ specializes to either $\ugl[k]$ for $\vpar\mapsto 1$ and to $\uqgl[k]$ for $\vpar\mapsto\qpar$.
The classical specialization is the one we studied in \autoref{SS:ProofClassical}.

We use the same notation as in \autoref{D:DualityQNum} 
for quantum numbers, but we see them as elements of $\aform$ 
or $\apara$ in general, 
and these specialize 
to the ones in \autoref{D:DualityQNum}.

The $\afrac$ algebra $\uvgl[k]$ is generated by
$E_{i}$, $F_{i}$ for $i\in\set{1,\dots,k-1}$,
and $L_{i}^{\pm 1}$ for $i\in\set{1,\dots,k}$ such that 
the $L_{i}^{\pm 1}$ commute with one another, 
$L_{i}^{-1}$ is the inverse of $L_{i}$, and
\begin{gather*}
L_{i}E_{j}=\vpar^{\delta_{i,j}-\delta_{i,j+1}}\cdot E_{j}L_{i}
,\quad
L_{i}F_{j}=\vpar^{-\delta_{i,j}+\delta_{i,j+1}}
\cdot F_{j}L_{i}
,\\
E_{i}F_{j}-F_{j}E_{i}=
\delta_{i,j}\tfrac{L_{i}L_{i+1}^{-1}-L_{i}^{-1}L_{i+1}}{\vpar-\vpar^{-1}}
,
\\
\vnum{2}\cdot E_{i}E_{j}E_{i}=E_{i}^{2}E_{j}+E_{j}E_{i}^{2}
,\;\text{if }|i-j|=1
,\\
E_{i}E_{j}-E_{j}E_{i}=0,\;\text{if }|i-j|>1
,
\\
\vnum{2}\cdot F_{i}F_{j}F_{i}=F_{i}^{2}F_{j}+F_{j}F_{i}^{2}
,\;\text{if }|i-j|=1
,\\
F_{i}F_{j}-F_{j}F_{i}=0,\;\text{if }|i-j|>1
,
\end{gather*}
for all suitable $i,j$.
We also choose the Hopf algebra structure on $\uagl[k]$ 
given by
\begin{alignat*}{3}
\Delta(E_{i})&=E_{i}\otimes L_{i}L_{i+1}^{-1}+1\otimes E_{i}
,\quad
\hspace{0.15cm}\epsilon(E_{i})&&=0
,\quad
S(E_{i})&&=-E_{i}L_{i}^{-1}L_{i+1},
\\
\Delta(F_{i})&=F_{i}\otimes 1+L_{i}^{-1}L_{i+1}\otimes F_{i}
,\quad
\epsilon(F_{i})&&=0
,\quad
S(F_{i})&&=-L_{i}L_{i+1}^{-1}F_{i},
\end{alignat*}
with $L_{i}^{\pm 1}$ being group like.

Following \cite{Lu-qgroups-root-of-1}, $\uagl[k]$ is the $\aform$
subalgebra of $\uvgl[k]$ generated by the \emph{divided powers} for $E_{i}$ and $F_{i}$, {\ie}
\begin{gather*}
E_{i}^{(j)}=\frac{E_{i}^{j}}{\vnum{j}!}
,
F_{i}^{(j)}=\frac{F_{i}^{j}}{\vnum{j}!}
,i\in\set{1,\dots,k-1},j\in\N,
\end{gather*}
and also by some adjustments of the $L_{i}$, see \cite{AnPoWe-representation-qalgebras}.
As the Hopf algebra structure of $\uagl[k]$ we take the one induced 
by $\uvgl[k]$.

\begin{Remark}\label{R:ProofQuantumSpecialize}
The $\aform$ algebra $\uagl[k]$ specializes to 
$\ugl[k]$ for $\vpar\mapsto 1$ and to 
$\uqgl[k]$ for $\vpar\mapsto\qpar$.
In both cases the divided power generators 
are only needed for $j=1$.
\end{Remark}

We scalar extend $\uagl[k]$ to an $\apara$ algebra, keeping the 
same notation.
The additional parameters only play a role for 
$\uagl[k]$ modules and not for $\uagl[k]$ itself.

\begin{Lemma}\label{L:ProofDualityVerma}
\autoref{D:DualityVerma} works {\ver} over $\apara$ 
giving $\uagl[2]$ modules.
\end{Lemma}

\begin{proof}
All appearing scalars can be interpreted in $\apara$.
\end{proof}

The $\uagl[2]$ modules from \autoref{L:ProofDualityVerma} are the 
\emph{integral Verma modules}. We denote these 
as before but using $\mathbb{A}$ as a subscript, {\eg} $\averma{\parav}$ 
is the integral version of $\qverma{\para}$. 
Using the Hopf algebra structure we can then define 
$\averma{\oplus\bsymv}$ similarly as we defined 
$\qverma{\oplus\bsym}$.

\begin{Remark}\label{R:ProofQuantumDerivations}
We now copy the approach taken in \autoref{SS:ProofClassical}. 
That is, we define a quantum polynomial algebra	$\apolyalg$
on which $\uagl[2]$ and $\uagl$ act by quantum 
derivatives in the spirit of {\eg} \cite[Section VII.3]{Ka-quantum-groups}.
This is done such that $\apolyalg\cong\averma{\oplus\bsymv}$ as $\uagl[2]$-$\uagl$ bimodules.
However, we do not use the language of quantum 
derivatives because 
of the various quantum parameters appear everywhere 
which make this setup cumbersome instead of helpful. For example, 
one would have relations of the form $Y_{j}X_{i}=X_{i}Y_{j}+(\vpar-\vpar^{-1})X_{j}Y_{i}$ 
and commutativity turns into quantum commutativity.
In order to avoid these technical difficulties we decided 
to define $\apolyalg$ instead as a free $\apara$ module with an explicit 
biaction defined on basis elements. The reader is still invited 
to think of the below as quantum 
derivatives acting on a quantum polynomial algebra.
\end{Remark}

\begin{Definition}\label{D:ProofQuantumDerivations}
We define the free $\apara$ module
\begin{gather*}
\apolyalg=\apara\{\overline{\bsym[X]}^{\bsymv+\bsym[r]}\overline{\bsym[Y]}^{\bsym[s]}|\bsym[r]\in\Z^{n},\bsym[s]\in\N^{n}\}
\end{gather*}
where we, as before, use formal parameters.
\end{Definition}

Write $X_{i}^{\parav_{i}+r_{i}}=
\vfac{\parav_{i}+r_{i}}\cdot\overline{X}_{i}^{\parav_{i}+r_{i}}$ 
and $Y_{i}^{s_{i}}=\vfac{s_{i}}\cdot\overline{Y}_{i}^{s_{i}}$.
We let $\uvgl[2]$ act on the scalar extension of $\apolyalg$ 
to $\afracpara$ by
\begin{gather}\label{Eq:ProofDiffOpQuantum}
\begin{aligned}
E\acts\bsym[X]^{\bsymv+\bsym[r]}\bsym[Y]^{\bsym[s]} 
&= 
\medsum_{i=1}^{n}
\vpar^{
\sums{\parav_{n}}+\sums{r_{n}}-\sums{\parav_{i}}-\sums{r_{i}}-\sums{s_{n}}+\sums{s_{i}}
}
\vnum{s_i}
\cdot
\bsym[X]^{\bsymv+\bsym[r]+\stdvec[i]}\bsym[Y]^{\bsym[s]-\stdvec[i]}
,
\\
F\acts\bsym[X]^{\bsymv+\bsym[r]}\bsym[Y]^{\bsym[s]} 
&= 
\medsum_{i=1}^{n}
\vpar^{-\sums{\parav_{i-1}}-\sums{r_{i-1}}+\sums{s_{i-1}}}
\vnum{\parav_{i}+r_{i}}
\cdot
\bsym[X]^{\bsymv+\bsym[r]-\stdvec[i]}\bsym[Y]^{\bsym[s]+\stdvec[i]}
,
\\
L_{1}\acts\bsym[X]^{\bsymv+\bsym[r]}\bsym[Y]^{\bsym[s]}
&=
\vpar^{\sums{\parav_{n}}+\sums{r_{n}}}
\cdot
\bsym[X]^{\bsymv+\bsym[r]}\bsym[Y]^{\bsym[s]}
,
\\
L_{2}\acts\bsym[X]^{\bsymv+\bsym[r]}\bsym[Y]^{\bsym[s]}
&=
\vpar^{\sums{s_{n}}}
\cdot
\bsym[X]^{\bsymv+\bsym[r]}\bsym[Y]^{\bsym[s]}
.
\end{aligned}
\end{gather}
We also define an $\uvgl[n]$ action 
on the scalar extension of $\apolyalg$ 
to $\afracpara$ by
\begin{gather}\label{Eq:ProofDiffOpQuantumTwo}
\begin{aligned}
\bsym[X]^{\bsymv+\bsym[r]}\bsym[Y]^{\bsym[s]}\acts E_{i}
=&
\vnum{\parav_{i+1}+r_{i+1}}\vpar^{\parav_{i}+r_{i}-\parav_{i+1}-r_{i+1}}
\cdot
\bsym[X]^{\bsymv+\bsym[r]+\sroot[i]}\bsym[Y]^{\bsym[s]}
\\&+
\vnum{s_{i+1}}
\cdot
\bsym[X]^{\bsymv+\bsym[r]}\bsym[Y]^{\bsym[s]+\sroot[i]}
,
\\
\bsym[X]^{\bsymv+\bsym[r]}\bsym[Y]^{\bsym[s]}\acts F_{i}
=&
\vnum{\parav_{i}+r_{i}}
\cdot
\bsym[X]^{\bsymv+\bsym[r]-\sroot[i]}\bsym[Y]^{\bsym[s]}
\\&+
\vpar^{-\parav_{i}-r_{i}+\parav_{i+1}+r_{i+1}}\vnum{s_{i}}
\cdot
\bsym[X]^{\bsymv+\bsym[r]}\bsym[Y]^{\bsym[s]-\sroot[i]}
,
\\
\bsym[X]^{\bsymv+\bsym[r]}\bsym[Y]^{\bsym[s]}\acts L_{i}
=&
\vpar^{\parav_{i}+r_{i}+s_{i}}\cdot\bsym[X]^{\bsymv+\bsym[r]}\bsym[Y]^{\bsym[s]}
.
\end{aligned}
\end{gather}

As before, we also get graded pieces $(\apolyalg)_{\bsym[d]}$ 
for $\bsym[d]\in\Z^{n}$.

\begin{Lemma}\label{L:ProofDiffOpQuantum}
The graded free $\aform$ module
\begin{gather*}
\apolyalg\cong\medoplus_{\bsym[d]\in\Z^{n}}(\apolyalg)_{\bsym[d]}
\end{gather*} 
is an $\uagl[2]$ module when endowed with 
the $\apara$-version of \autoref{Eq:ProofDiffOpQuantum}
that decomposes as above and 
is isomorphic to $\averma{\oplus\bsymv}$.
Moreover, it is also an $\uagl$ module when endowed 
with the $\apara$-version of \autoref{Eq:ProofDiffOpQuantumTwo},
and the two actions commute.
\end{Lemma}

\begin{proof}
A direct calculation verifies that 
\autoref{Eq:ProofDiffOpQuantum} and \autoref{Eq:ProofDiffOpQuantumTwo}
give the scalar extension of $\apolyalg$ the structure of an
$\uvgl[2]$-$\uvgl$ bimodule. 
One then checks that \autoref{Eq:ProofDiffOpQuantum} and \autoref{Eq:ProofDiffOpQuantumTwo} when successively applied 
to $E_{i}^{(j)}$ and 
$F_{i}^{(j)}$ has coefficients in $\apara$ 
when acting on the basis 
of the elements of the form $\overline{\bsym[X]}^{\bsymv+\bsym[r]}\overline{\bsym[Y]}^{\bsym[s]}$ 
(to see this it is helpful to keep the quantum derivative 
picture from \autoref{R:ProofQuantumDerivations} and 
$\partial_{X}\frac{X^{r}}{r!}=\frac{X^{r-1}}{(r-1)!}$
in mind), so the 
biaction can be restricted and gives an $\uagl[2]$-$\uagl$ bimodule 
structure on $\apolyalg$.

The final claim that 
$\apolyalg\cong\averma{\oplus\bsymv}$ as $\uagl[2]$-$\uagl$ bimodules 
can be verified as in \autoref{L:ProofDiffOp}.
\end{proof}

\begin{Definition}\label{D:ProofFlat}
\leavevmode
\begin{enumerate}

\item We call an $\apara$ module $\module[M]$ \emph{(generically) flat} 
if its free and 
all specializations to characteristic zero fields where $\bsymv$ 
are specialized to admissible parameters are of the same dimension.

\item We call an $\uagl[k]$ module $\module[M]$ 
\emph{(generically) flat} if it is flat as an $\apara$ module and 
if all specializations of 
$\End_{\uagl[k]}(\module[M])$ to characteristic zero fields where $\bsymv$ 
are specialized to admissible parameters are of the same dimension.

\item By being \emph{(generically) flat as an $\uagl[2]$-$\uagl$ bimodule} 
we mean being flat as an $\uagl[2]$ module and as an $\uagl$ module.

\end{enumerate}
\end{Definition}

\begin{Lemma}\label{L:ProofFlat}
The $\uagl[2]$-$\uagl$ bimodule 
$\averma{\oplus\bsymv}\cong\apolyalg$ is flat, 
its classical specialization is 
the $\ugl[2]$-$\ugl$ bimodule 
$\verma{\oplus\bsym}\cong\polyalg$
and its quantum specialization is 
the $\uqgl[2]$-$\uqgl$ bimodule 
$\qverma{\oplus\bsym}\cong\qpolyalg$.
\end{Lemma}

\begin{proof}
The quantum version of \autoref{L:ProofSemisimple} 
holds as well, with the same proof. This implies that 
$\averma{\bsymv}$ is tilting when specialized to characteristic zero fields 
with $\bsymv\mapsto\bsym$ for $\bsym$ admissible parameters. Therefore $\averma{\bsymv}$ is flat 
by the usual arguments, see {\eg} the 
arXiv appendix to \cite{AnStTu-cellular-tilting}.
Moreover, comparison of formulas implies that 
the specializations are the claimed ones.
\end{proof}

\begin{Remark}\label{R:ProofTypeI}
For the below note that, 
by their very construction, all 
$\uagl[k]$ modules used in this paper are of \emph{type 
$(1,\dots,1)$} in the sense of {\eg} 
\cite[Section 1.4]{AnPoWe-representation-qalgebras} or
\cite[Section 5.2]{Ja-lectures-qgroups}.
\end{Remark}

\begin{proof}[Proof of \autoref{T:DualityHowe}]
As we will see now, 
using flatness, the quantum Verma Howe duality theorem
follows from the classical case.
Our exposition below follows \cite[Section 7A]{SaTu-bcd-webs},
but flatness arguments 
along the same lines are very common in the literature.	

\textit{Commuting actions.}
To use \autoref{L:ProofFlat}, one first needs to establish the 
existence of the commuting actions as 
in \autoref{T:DualityHowe} in the quantum 
case independently of the classical 
case. This is done in \autoref{L:ProofDiffOpQuantum}, so we can focus 
on the $\uqgl[2]$-$\uqgl$ bimodule decomposition.

As before in the 
classical case, all parameters in this proof are admissible from now on.

\textit{Bimodule decomposition.}
We will now repeatedly use \autoref{L:ProofFlat}.
We compare the $\uqgl[2]$ module $\qverma{\oplus\bsym}$
and the $\ugl[2]$ module $\verma{\oplus\bsym}$, and we see that 
the weights of these modules are the same under the usual identification
of quantum and classical weights. Moreover, the weight multiplicities 
are also the same and all finite. It follows 
then from \autoref{L:ProofFlat} that we have
\begin{gather*}
\qverma{\oplus\bsym}\cong
\medoplus_{\substack{\countingpara\in\Z\\ t\in\N}}
\qverma{\sums{\para_{n}}+\countingpara-t,t}\otimes\qdensemult{\sums{\para_{n}}+\countingpara-t,t}
\end{gather*}
as the quantum analog of \autoref{Eq:ProofDecomTensor}, where 
$\qdensemult{\sums{\para_{n}}+\countingpara-t,t}$ are 
multiplicity $\Kpara$ vector spaces. 
We actually know that these multiplicity 
$\Kpara$ vector spaces are $\uqgl$ modules by the 
quantum specialization of the
previously established $\uvgl[2]$-$\uvgl$ bimodule 
structure.

We want to show that all appearing 
$\qdensemult{\sums{\para_{n}}+\countingpara-t,t}$ are simple
as $\uqgl$ modules. This is equivalent to
the action giving a surjection
\begin{gather}\label{Eq:ProofSurjection}
f\colon\uqgl\twoheadrightarrow
\End_{\uqgl[2]}(\qverma{\bsym}).
\end{gather}
Now, setting $\vpar\mapsto1$ or $\vpar\mapsto\qpar$
, respectively, and $\bsymv\mapsto\bsym$ in the $\apara$ 
version, we can identify 
$\averma{\oplus\bsymv}$ with 
$\verma{\oplus\bsym}$, and the biactions of
$\uagl[2]/(\vpar-1,\bsymv-\bsym)$-$\uagl/(\vpar-1,\bsymv-\bsym)$ 
and $\ugl[2]$-$\ugl$ coincide under this specialization, and 
{\ver} for $\uqgl[2]$-$\uqgl$ instead of 
$\ugl[2]$-$\ugl$.
In particular, the images of 
these two actions agree. It follows now from 
the classical version of \autoref{T:DualityHowe} that the 
action map $f$ is surjective classically. 
Thus, \autoref{L:ProofFlat} implies that 
\autoref{Eq:ProofSurjection} holds and 
$\qdensemult{\sums{\para_{n}}+\countingpara-t,t}$ are simple
as $\uqgl$ modules.

Comparison of definitions verifies that the classical version of $\qdensemult{\sums{\para_{n}}+\countingpara-t,t}$ 
and $\qdense{\countingpara-t,t}$ are $\densemult{\sums{\para_{n}}+\countingpara-t,t}$ 
and $\dense{\countingpara-t,t}$, respectively. By the classification recalled in
\cite[Section 2.3]{Ma-gt-dense} (originally proven in \cite{Ma-classification-weight-modules})
we get also that $\densemult{\sums{\para_{n}}+\countingpara-t,t}\cong\dense{\countingpara-t,t}$, 
as before.
Finally, $\qdensemult{\sums{\para_{n}}+\countingpara-t,t}$ 
and $\qdense{\countingpara-t,t}$ are 
of type $(1,\dots,1)$, by construction,
and true quantum deformation in the sense of  
\cite[Section 2.1]{Ma-gt-dense} which implies 
$\qdensemult{\sums{\para_{n}}+\countingpara-t,t}\cong\qdense{\countingpara-t,t}$.

\textit{Dense.} 
By \autoref{L:ProofFlat} 
and the $\uqgl[2]$-$\uqgl[n]$ bimodule decomposition 
from \autoref{Eq:DualityDecom}, the argument is now the same as in the 
classical case.
\end{proof}

%%%%%%%%%%%%%%%%%%%%%%%%%%%%%%%%%%%%%%%%%

\section{The colored higher LKB representations are simple}\label{S:LKB}

%%%%%%%%%%%%%%%%%%%%%%%%%%%%%%%%%%%%%%%%%

Recall that we have fixed $n\in\Z_{\geq 1}$ and parameters $\bsym=(\para_{1},\dots,\para_{n})\in\K^{n}$.

%%%%%%%%%%%%%%%%%%%%%%%%%%%%%%%%%%%%%%%%%

\subsection{Pure and colored braids}\label{SS:LKBBraids}

%%%%%%%%%%%%%%%%%%%%%%%%%%%%%%%%%%%%%%%%%

Let $\braid$ denote the \emph{braid group} with $n$ strands
which can be illustrated using the usual diagrammatics, {\eg}
\begin{gather}\label{Eq:LKBBraidExample}
\begin{tikzpicture}[anchorbase,scale=0.55]
\draw[usual,crossline,rounded corners] (7,0) to[out=90,in=0] (6.5,0.2) to (4.5,0.2) to (4.5,0.8) to[out=90,in=270] (4,2);
\draw[usual,crossline] (1,0) to[out=90,in=270] (1,1);
\draw[usual,crossline] (0,0) to[out=90,in=270] (4,1);
\draw[usual,crossline] (2,0) to[out=90,in=270] (3,2);
\draw[usual,crossline] (4,1) to[out=90,in=270] (0,2);
\draw[usual,crossline] (6,1) to[out=90,in=270] (5,2);
\draw[usual,crossline] (6,0) to[out=90,in=270] (5,1) to[out=90,in=270] (6,2);
\draw[usual,crossline] (3,0) to[out=90,in=180] (4,1) to[out=0,in=270] (8,2);
\draw[usual,crossline] (4,0) to[out=90,in=270] (2,2);
\draw[usual,crossline] (5,0) to[out=90,in=270] (6,1);
\draw[usual] (6,0.95) to (6,1.05);
\draw[usual,crossline] (8,0) to[out=90,in=270] (7,2);
\draw[usual,crossline] (1,1) to[out=90,in=270] (1,2);
\end{tikzpicture}
\,,\quad
\bgen[{i,i+1}]=
\begin{tikzpicture}[anchorbase,scale=0.55]
\draw[usual] (1,0)node[below]{$i{+}1$} to[out=90,in=270] (0,1)node[above,yshift=-1pt]{$i{+}1$};
\draw[usual,crossline] (0,0)node[below]{$i$} to[out=90,in=270] (1,1)node[above,yshift=-1pt]{$i$};
\end{tikzpicture}
\,,\quad
\bgen[{i,i+1}]^{-1}=
\begin{tikzpicture}[anchorbase,scale=0.55]
\draw[usual] (0,0)node[below]{$i$} to[out=90,in=270] (1,1)node[above,yshift=-1pt]{$i$};
\draw[usual,crossline] (1,0)node[below]{$i{+}1$} to[out=90,in=270] (0,1)node[above,yshift=-1pt]{$i{+}1$};
\end{tikzpicture}
.
\end{gather}
As displayed above, the transposition generators, 
crossing the $i$th and $(i+1)$th strand, of $\braid$ 
are denoted by $\bgen[{i,i+1}]$ and $\bgen[{i,i+1}]^{-1}$.

Recall that the \emph{pure braid group} is the subgroup $\pbraid\subset\braid$
of all elements with the bottom and the top of each strand in the same position. More generally, we define:

\begin{Definition}\label{D:LKBPureBraids}
Let $\setstuff{P}(\{1,\dots,n\})$ be the set of partitions of $\{1,\dots,n\}$. For every $S\in\setstuff{P}(\{1,\dots,n\})$, 
the \emph{braid group} that is \emph{pure on $S$} is the subgroup $\ibraid{S}\subset\braid$
such that the strands with bottom points in $A\in S$ 
have their top points in $A$ as well.
\end{Definition}

\begin{Example}\label{E:LKBBraids}
We have $\ibraid{\{[1],\dots,[n]\}}=\pbraid$ and 
$\ibraid{\{[1,\dots,n]\}}=\braid$, where we use square 
brackets for the parts of the partition.	
Moreover, the leftmost braid in
\autoref{Eq:LKBBraidExample} is pure on the partition
$S=\{[1],[2],[3,4,5,8,9],[6],[7]\}$, and $S$ is the finest 
partition such that the braid is pure on it.
\end{Example}

\begin{Example}\label{E:LKBBraidsHandlebody}
The \emph{handlebody braid group} of genus $g\in\N$ with $n\in\Z_{\geq 1}$ strands is the subgroup of $\braid[{g+n}]$ that is pure on
$S=\{[1],\dots,[g],[g+1,\dots,g+n]\}$. (For $g=0$, by convention, $\braid[{g+n}]$ is the classical braid group $\braid[{n}]$.) The first $g$ strands in the handlebody braid group are core strands, while the remaining strands are usual strands, see e.g. \cite[Section 2]{RoTu-homflypt-handlebody} for the topological background.
\end{Example}

\begin{Definition}\label{D:LKBBraids}
We associate a partition $S(\bsym)\in\setstuff{P}(\{1,\dots,n\})$ to 
$\bsym$ by 
\begin{gather*}
\text{$i,j$ are in the same component of $S(\bsym)$ $\Leftrightarrow$ 
$\para_{i}=\para_{j}$.}
\end{gather*}
We denote the 
corresponding subgroup of the braid group 
by $\ibraid{\bsym}=\ibraid{S(\bsym)}$.
\end{Definition}

\begin{Example}\label{E:LKBBraidsTwo}
If all $\para_{i}$ are different, then 
$\ibraid{\bsym}=\pbraid$, and if 
$\bsym=(\para,\dots,\para)$, then 
$\ibraid{\bsym}=\braid$. For the leftmost braid in
\autoref{Eq:LKBBraidExample} the finest set of parameters 
is $\bsym=(\para_{1},\para_{2},\para_{3},
\para_{3},\para_{3},\para_{6},\para_{7},\para_{3},\para_{3})$ 
for pairwise distinct $\para_{i}$.
\end{Example}

\begin{Example}\label{E:LKBBraidsHandlebodyTwo}
For the handlebody braid group as in \autoref{E:LKBBraidsHandlebody}, a natural choice of $\bsym$ is $\para_{1}=\dots=\para_{g}\notin\Z$ and 
$\para_{g+1}=\dots=\para_{g+n}\in\Z$ otherwise. Note that such a 
choice of parameters is admissible 
for $g>0$.
\end{Example}

%%%%%%%%%%%%%%%%%%%%%%%%%%%%%%%%%%%%%%%%%

\subsection{LKB representations}\label{SS:LKBSimple}

%%%%%%%%%%%%%%%%%%%%%%%%%%%%%%%%%%%%%%%%%

We again work over $\apara$ and $\Kpara$. We always consider 
$\averma{\bsymv}=\averma{\parav_{1}}
\otimes\dots\otimes\averma{\parav_{n}}$ as a $\uagl[2]$ module via \autoref{Eq:ProofDiffOpQuantum} (note that $\apolyalg\cong\averma{\oplus\bsymv}
=\medoplus_{\bsym[d]\in\Z^{n}}\averma{\parav_{1}+d_{1}}
\otimes\dots\otimes\averma{\parav_{n}+d_{n}}$, but the direct sum is
only needed for the dual action).

We now adjust the construction from 
\cite[Section 3]{JaKe-verma-lkb} 
(and the references to \cite{Ka-quantum-groups} therein)
to our setting. See also \cite[Definition 2.19]{Ma-colored-lkb}.

\begin{Definition}\label{D:LKBRMatrix}
Let $\lkbscale[\pm 1]\in\End_{\apara}(\averma{\parav_{i}}\otimes
\averma{\parav_{j}})$ defined by
\begin{gather*}
\lkbscale[\pm 1](m_{k}\otimes m_{l})
=
\vpar^{\pm(-l\parav_{i}-k\parav_{j}+2kl)}
\cdot(m_{k}\otimes m_{l}).
\end{gather*}
Write $F^{[r]}=\frac{(\vpar-\vpar^{-1})^{r}}{\vfac{r}}F^{r}$ 
for $r\in\N$.
Define the \emph{R matrix} and its inverse on 
$\averma{\parav_{1}}\otimes
\averma{\parav_{2}}$ as
\begin{align*}
\rmatrix[{\parav_{i},\parav_{j}}]
&\colon\averma{\parav_{i}}\otimes
\averma{\parav_{j}}\to\averma{\parav_{j}}\otimes
\averma{\parav_{i}},
\quad
\rmatrix[{\parav_{i},\parav_{j}}]=
s\circ\lkbscale[+1]\circ\bigg(\medsum_{l=0}^{\infty}
\vpar^{l(l-1)/2}\cdot
e^{l}\otimes f^{[l]}\bigg),
\\
\rimatrix[{\parav_{i},\parav_{j}}]
&\colon\averma{\parav_{i}}\otimes
\averma{\parav_{j}}\to\averma{\parav_{j}}\otimes
\averma{\parav_{i}},
\quad
\scalebox{0.8}{$\rimatrix[{\parav_{i},\parav_{j}}]=
\bigg(\medsum_{l=0}^{\infty}
(-1)^{l}\vpar^{-l(l-1)/2}\cdot
e^{l}\otimes f^{[l]}\bigg)\circ\lkbscale[-1]\circ s,$}
\end{align*}
where $s$ is the swap map $s(x\otimes y)=y\otimes x$.
\end{Definition}

\begin{Lemma}\label{L:LKBFiniteOperator}
The operators $\rmatrix[{\parav_{i},\parav_{j}}]$ and 
$\rimatrix[{\parav_{i},\parav_{j}}]$ are well-defined, {\ie} 
the appearing summations are finite on every $m_{k}\otimes m_{l}$. 
\end{Lemma}

\begin{proof}
This holds because the operator $e$ is locally nilpotent.
\end{proof}

Graphically we will denote these operators by
\begin{gather}\label{Eq:LKBRMatrix}
\rmatrix[{\parav_{i},\parav_{j}}]
\leftrightsquigarrow
\begin{tikzpicture}[anchorbase,scale=0.55]
\draw[usual] (1,0)node[below]{$\parav_{j}$} to[out=90,in=270] (0,1)node[above,yshift=-2.7pt]{$\parav_{j}$};
\draw[usual,crossline] (0,0)node[below]{$\parav_{i}$} to[out=90,in=270] (1,1)node[above,yshift=-1.1pt]{$\parav_{i}$};
\end{tikzpicture}
,\quad
\rimatrix[{\parav_{i},\parav_{j}}]
\leftrightsquigarrow
\begin{tikzpicture}[anchorbase,scale=0.55]
\draw[usual] (0,0)node[below]{$\parav_{i}$} to[out=90,in=270] (1,1)node[above,yshift=-1.1pt]{$\parav_{i}$};
\draw[usual,crossline] (1,0)node[below]{$\parav_{j}$} to[out=90,in=270] (0,1)node[above,yshift=-2.7pt]{$\parav_{j}$};
\end{tikzpicture}
.
\end{gather}
We now define a $\ibraid{\bsym}$ action on 
$\averma{\bsym}$ by \emph{colored reading}. 
That is, one colors the strands 
of $\beta\in\ibraid{\bsym}$ by $\bsym$, and then 
we get an element of $\End_{\apara}(\averma{\bsym})$ 
by composing the relevant version of 
\autoref{Eq:LKBRMatrix} from bottom to top.
We call this element $\rmatrix[\beta]$.

\begin{Example}\label{E:LKBRMatrix}
For $\bsym=(\parav_{1},\parav_{2},\parav_{3},
\parav_{3},\parav_{3},\parav_{6},\parav_{7},\parav_{3},\parav_{3})$
and the leftmost braid in
\autoref{Eq:LKBBraidExample}
we get
\begin{gather*}
\begin{tikzpicture}[anchorbase,scale=0.55]
\draw[usual,crossline,rounded corners] (7,0)node[below]{$\parav_{3}$} to[out=90,in=0] (6.5,0.2) to (4.5,0.2) to (4.5,0.8) to[out=90,in=270] (4,2)node[above,yshift=-1pt]{$\parav_{3}$};
\draw[usual,crossline] (1,0)node[below]{$\parav_{2}$} to[out=90,in=270] (1,1);
\draw[usual,crossline] (0,0)node[below]{$\parav_{1}$} to[out=90,in=270] (4,1);
\draw[usual,crossline] (2,0)node[below]{$\parav_{3}$} to[out=90,in=270] (3,2)node[above,yshift=-1pt]{$\parav_{3}$};
\draw[usual,crossline] (4,1) to[out=90,in=270] (0,2)node[above,yshift=-1pt]{$\parav_{1}$};
\draw[usual,crossline] (6,1) to[out=90,in=270] (5,2)node[above,yshift=-1pt]{$\parav_{6}$};
\draw[usual,crossline] (6,0)node[below]{$\parav_{7}$} to[out=90,in=270] (5,1) to[out=90,in=270] (6,2)node[above,yshift=-1pt]{$\parav_{7}$};
\draw[usual,crossline] (3,0)node[below]{$\parav_{3}$} to[out=90,in=180] (4,1) to[out=0,in=270] (8,2)node[above,yshift=-1pt]{$\parav_{3}$};
\draw[usual,crossline] (4,0)node[below]{$\parav_{3}$} to[out=90,in=270] (2,2)node[above,yshift=-1pt]{$\parav_{3}$};
\draw[usual,crossline] (5,0)node[below]{$\parav_{6}$} to[out=90,in=270] (6,1);
\draw[usual] (6,0.95) to (6,1.05);
\draw[usual,crossline] (8,0)node[below]{$\parav_{3}$} to[out=90,in=270] (7,2)node[above,yshift=-1pt]{$\parav_{3}$};
\draw[usual,crossline] (1,1) to[out=90,in=270] (1,2)node[above,yshift=-1pt]{$\para_{2}$};
\end{tikzpicture}
\\
\rightsquigarrow
\rmatrix[\beta]=
\rmatrix[{\parav_{2},\parav_{1}}]\circ\rmatrix[{\parav_{3},\parav_{1}}]
\circ\dots\circ
\rimatrix[{\parav_{1},\parav_{3}}]\circ\rmatrix[{\parav_{1},\parav_{2}}]
\in\End_{\Kpara}(\averma{\bsym})
.
\end{gather*}
The endomorphism $\rmatrix[\beta]$ has 
eighteen $R$ matrix factors in total.
\end{Example}

\begin{Definition}\label{D:LKBRefinement}
A \emph{refinement} $\bsym[\rho]$ of $\bsymv$ is a set 
of parameters that gives a refined partition compared to $\bsymv$
when applying \autoref{D:LKBBraids}.
We write $\bsym[\rho]\leq\bsymv$ for refinements of $\bsymv$.
\end{Definition}

\begin{Notation}\label{N:LKBAction}
For $\uvgl[2]$ we extend scalars to 
$\afrac$ or $\afracpara$ but do not indicate 
this in the notation.
\end{Notation}

Similar to the braid group action in symmetric Howe duality, 
the braid group acts on one $\gl[n]$ weight space of 
$\averma{\oplus\bsym}$ and this action commutes with the 
$\gl[2]$ action:

\begin{Lemma}\label{L:LKBAction}
\leavevmode

\begin{enumerate}

\item \autoref{Eq:LKBRMatrix} and colored reading 
endows $\averma{\bsymv}$ with the structure 
of a $\ibraid{\bsym[\rho]}$ module for $\bsym[\rho]\leq\bsymv$.

\item \scalebox{0.95}{Colored reading commutes with the $\uvgl[2]$ action 
coming from \autoref{Eq:ProofDiffOpQuantum}.}

\item The image of $\ibraid{\bsym[\rho]}$ under this module 
structure is in $\End_{\uvgl[2]}(\averma{\bsymv})$.

\end{enumerate}
\end{Lemma}

\begin{proof}
One first proves, by copying 
\cite[Theorem 7]{JaKe-verma-lkb}, 
that \autoref{Eq:LKBRMatrix} satisfies the colored 
braid relations, {\eg}
\begin{gather*}
\begin{tikzpicture}[anchorbase,scale=0.55]
\draw[usual] (1,0)node[below]{$\parav_{j}$} to[out=90,in=270] (0,1) to[out=90,in=270] (1,2)node[above,yshift=-2.7pt]{$\parav_{j}$};
\draw[usual,crossline] (0,0)node[below]{$\parav_{i}$} to[out=90,in=270] (1,1) to[out=90,in=270] (0,2)node[above,yshift=-1.1pt]{$\parav_{i}$};
\end{tikzpicture}
=
\begin{tikzpicture}[anchorbase,scale=0.55]
\draw[usual] (1,0)node[below]{$\parav_{j}$} to[out=90,in=270] (1,2)node[above,yshift=-2.7pt]{$\parav_{j}$};
\draw[usual,crossline] (0,0)node[below]{$\parav_{i}$} to[out=90,in=270] (0,2)node[above,yshift=-1.1pt]{$\parav_{i}$};
\end{tikzpicture}
,\quad
\begin{tikzpicture}[anchorbase,scale=0.55]
\draw[usual] (1,0)node[below]{$\parav_{j}$} to[out=90,in=270] (0,1) to[out=90,in=270] (0,2)node[above,yshift=-2.7pt]{$\parav_{j}$};
\draw[usual,crossline] (0,0)node[below]{$\parav_{i}$} to[out=90,in=270] (1,1) to[out=90,in=270] (1,2)node[above,yshift=-1.1pt]{$\parav_{i}$};
\draw[usual] (3,0)node[below]{$\parav_{l}$} to[out=90,in=270] (3,1) to[out=90,in=270] (2,2)node[above,yshift=-1.25pt]{$\parav_{l}$};
\draw[usual,crossline] (2,0)node[below]{$\parav_{k}$} to[out=90,in=270] (2,1) to[out=90,in=270] (3,2)node[above,yshift=-1.25pt]{$\parav_{k}$};
\end{tikzpicture}
=
\begin{tikzpicture}[anchorbase,scale=0.55]
\draw[usual] (1,0)node[below]{$\parav_{j}$} to[out=90,in=270] (1,1) to[out=90,in=270] (0,2)node[above,yshift=-2.7pt]{$\parav_{j}$};
\draw[usual,crossline] (0,0)node[below]{$\parav_{i}$} to[out=90,in=270] (0,1) to[out=90,in=270] (1,2)node[above,yshift=-1.1pt]{$\parav_{i}$};
\draw[usual] (3,0)node[below]{$\parav_{l}$} to[out=90,in=270] (2,1) to[out=90,in=270] (2,2)node[above,yshift=-1.25pt]{$\parav_{l}$};
\draw[usual,crossline] (2,0)node[below]{$\parav_{k}$} to[out=90,in=270] (3,1) to[out=90,in=270] (3,2)node[above,yshift=-1.25pt]{$\parav_{k}$};
\end{tikzpicture}
,\\
\begin{tikzpicture}[anchorbase,scale=0.55]
\draw[usual] (2,0)node[below]{$\parav_{k}$} to[out=90,in=270] (0,2)node[above,yshift=-1.25pt]{$\parav_{k}$};
\draw[usual,crossline] (1,0)node[below]{$\parav_{j}$} to[out=90,in=270] (0,1) to[out=90,in=270] (1,2)node[above,yshift=-2.7pt]{$\parav_{j}$};
\draw[usual,crossline] (0,0)node[below]{$\parav_{i}$} to[out=90,in=270] (2,2)node[above,yshift=-1.1pt]{$\parav_{i}$};
\end{tikzpicture}
=
\begin{tikzpicture}[anchorbase,scale=0.55]
\draw[usual] (2,0)node[below]{$\parav_{k}$} to[out=90,in=270] (0,2)node[above,yshift=-1.25pt]{$\parav_{k}$};
\draw[usual,crossline] (1,0)node[below]{$\parav_{j}$} to[out=90,in=270] (2,1) to[out=90,in=270] (1,2)node[above,yshift=-2.7pt]{$\parav_{j}$};
\draw[usual,crossline] (0,0)node[below]{$\parav_{i}$} to[out=90,in=270] (2,2)node[above,yshift=-1.1pt]{$\parav_{i}$};
\end{tikzpicture}
.
\end{gather*}
Hence, $\rmatrix[\beta]$ is independent
of the choices in colored reading, 
and we obtain the claimed $\ibraid{\bsym[\mu]}$ module
structure. For
$\ibraid{\bsym[\rho]}\subset\ibraid{\bsymv}$ 
this $\ibraid{\bsymv}$ module structure 
restricts to $\ibraid{\bsym[\rho]}$.

That the two actions commute follows because of the well-known fact 
(and easy calculation) that the R matrices 
are $\uvgl[2]$ equivariant with respect to \autoref{Eq:ProofDiffOpQuantum}.

The final claim follows since the action maps 
commute with the $\uvgl[2]$ action on $\averma{\bsymv}$.
\end{proof}

We thus have a $\uvgl[2]$-$\ibraid{\bsym[\rho]}$ bimodule structure on $\averma{\bsymv}$.

\begin{Remark}\label{R:LKBHalfQG}
\autoref{L:LKBAction} can be strengthened: the image 
of $\ibraid{\bsym[\rho]}$ commutes with the action of 
the $\apara$ subalgebra of $\uagl[2]$ generated by 
$E$, $F^{[r]}$ for $r\in\N$, $L_{1}^{\pm 1}$ and $L_{2}^{\pm 1}$.
\end{Remark}

We now turn our attention to the 
(colored higher) LKB representations, which, 
following \cite[Section 3]{JaKe-verma-lkb}, we define as follows:

\begin{Definition}\label{D:LKBReps}
Let $\ker(e)$ and $\ker(k-\medprod_{i=1}^{n}\vpar^{\parav_{i}}\vpar^{-2l})
=\ker(k-\vpar^{\sums{\parav_{n}}-2l})$ 
be the kernels of the indicated operators coming from 
the $\uvgl[2]$ action on $\averma{\bsymv}$.
For $l\in\N$ the \emph{$l$th LKB representation} is
defined as $\alkb{{n,l}}=\ker(e)\cap
\ker(k-\vpar^{\sums{\parav_{n}}-2l})$.
\end{Definition}

\begin{Lemma}\label{L:LKBReps}
For $\bsym[\rho]\leq\bsymv$ the free $\apara$ module
$\alkb{{n,l}}$ is stable under the $\ibraid{\bsym[\rho]}$ action.
\end{Lemma}

\begin{proof}
Note first that the definition of $\alkb{{n,l}}$ only involves the 
operators $e$, $k$ and $\apara$ multiples of the identity. 
Hence, $\alkb{{n,l}}$ is defined over $\apara$, by construction.	
The rest can be proven, {\muta}, as in \cite[Theorem 1]{JaKe-verma-lkb}.
\end{proof}

\begin{Example}\label{E:LKBRepsHW}
With respect to the basis $\{m_{i}|i\in\N\}$ in \autoref{E:DualityVerma},
$\alkb{{n,1}}$ has a basis given by an $n$ fold tensor product of 
$m_{i}$ with one entry being $m_{1}$ and all other entries being 
$m_{0}$. Hence, the $\apara$ rank 
is $\binom{n-1}{1}$. In general, $\alkb{{n,l}}$ 
is of $\apara$ rank $\binom{n+l-2}{l}$.
\end{Example}

\begin{Remark}\label{R:LKBParameters}
The representation $\alkb{{n,l}}$ has between two and 
$n+1$ parameters, depending on $\bsymv$. For 
example, for $\bsymv=(\parav,\dots,\parav)$ 
one has $\vpar$ and $\parav$ as parameters.
\end{Remark}

\begin{Example}\label{E:LKBReps}
The representation $\alkb{{n,0}}$ is always trivial, while $\alkb{{n,1}}$ 
is the (reduced) \emph{Burau representation} of 
$\braid$, and $\alkb{{n,2}}$ is its classical 
\emph{LKB representation} as in \cite{La-lkb-rep}, 
or closer to our formulation, as in \cite{JaKe-verma-lkb}.

Or, to be completely precise, $\alkb{{n,2}}$ is
a multiparameter version of the construction from
\cite{JaKe-verma-lkb}, see also \cite{Ma-colored-lkb}. 
Moreover, the representation $\alkb{{n,l}}$ can then further be matched with 
it homological counterpart 
up to playing with parameters, see \cite[Theorem 6.1]{Ko-quantum-braid-lkb} 
and \cite[Theorem 1.5]{Ma-homological-colored-lkb}
for a precise statement.

Using an appropriate ground field 
and quantum parameter,
$\qlkb{{n,2}}$ for $\bsym=(\para,\dots,\para)$ is 
a faithful $\braid$ module by \cite{Bi-linear-artin} and 
\cite{Kr-linear-artin}, and thus, $\qlkb{{n,2}}$ is also 
faithful for $\ibraid{\bsym[\rho]}$ for all $\bsym[\rho]$.
\end{Example}

Specializing to the quantum case, the following is our 
main application of \autoref{T:DualityHowe}:

\begin{Theorem}\label{T:LKBSimple}
Assume that the parameters are admissible.
Then the representation $\qlkb{{n,l}}$ is 
a simple $\ibraid{\bsym[\rho]}$ module 
for $\bsym[\rho]\leq\bsym$ and all $l\in\N$.
\end{Theorem}

The proof of \autoref{T:LKBSimple}
is given in \autoref{S:ProofSimple}.

\begin{Remark}\label{R:LKBSimple}
\autoref{T:LKBSimple} extends and 
generalizes \cite[Theorem 3]{JaKe-verma-lkb} in multiple ways. First, \autoref{T:LKBSimple} is a multiparameter
version of \cite[Theorem 3]{JaKe-verma-lkb}. \autoref{T:LKBSimple} also generalizes the result in {\loccit} to arbitrary fields and generic $\qpar$, and also allows much more general parameters.
And even when 
we have only one parameter, {\ie} $\bsym=(\para,\dots,\para)$, and 
work over $\Q(\para,q)$ 
\autoref{T:LKBSimple} is stronger 
than \cite[Theorem 3]{JaKe-verma-lkb} 
since we {\eg} also prove that $\qlkb{{n,l}}$ is 
a simple $\pbraid$ module not just a simple $\braid$ module.
The proof given below is also very different from the 
one given in \cite{JaKe-verma-lkb}, and 
we do not know how to generalize the proof 
in \cite{JaKe-verma-lkb} 
to {\eg} include the various subgroups of $\braid$, including the pure and handlebody braid groups.
\end{Remark}

%%%%%%%%%%%%%%%%%%%%%%%%%%%%%%%%%%%%%%%%%

\section{The proof of simplicity}\label{S:ProofSimple}

%%%%%%%%%%%%%%%%%%%%%%%%%%%%%%%%%%%%%%%%%

Our proof of \autoref{T:LKBSimple} 
uses Verma versions of \cite[Theorem 5.5 and Remark 8.6]{LeZh-strongly-multiplicity-free}.

\begin{Remark}\label{R:ProofSimpleSMF}
We think of Verma modules as 
limits of symmetric powers and this
was one of our main motivation to follow the approach
taken in \cite{LeZh-strongly-multiplicity-free}.
\end{Remark}

%%%%%%%%%%%%%%%%%%%%%%%%%%%%%%%%%%%%%%%%%

\subsection{The classical case}\label{SS:ProofSimpleClassical}

%%%%%%%%%%%%%%%%%%%%%%%%%%%%%%%%%%%%%%%%%

Our ground field 
in this section is $\K$.

\begin{Definition}\label{D:ProofSimpleInfBraids}
For $k\in\Z_{\geq 2}$ 
the \emph{infinitesimal pure braid group} is the $\K$ algebra $\infbraid$ 
generated by $\infbgen$ for $1\leq i<j\leq k$ subject to
\begin{gather*}
[\infbgen,\infbgen[{rs}]]
=
[\infbgen[{ir}]+\infbgen[{is}],\infbgen[{rs}]]
=
[\infbgen,\infbgen[{ir}]+\infbgen[{jr}]]
=0,
\end{gather*}
for pairwise distinct $i,j,r,s$, and $[\placeholder,\placeholder]$ 
denotes the commutator. For $k=1$ we let $\infbraid=\K$.
\end{Definition}

\begin{Remark}\label{R:ProofSimpleInfBraids}
The motivation to study $\infbraid$ is that it gives 
rise to the so-called \emph{monodromy representation 
of the KZ equation} of the pure braid group $\pbraid$ 
which works for very general tensor products of Lie algebra 
representations, see 
\cite[Proposition 2.3]{Ko-cft-topology} for details.
\end{Remark}

\begin{Definition}\label{D:ProofSimpleOperators}
For all $1\leq i<j\leq n$ define operators on 
$\verma{\oplus\bsym}\cong\polyalg$ by
\begin{gather}\label{Eq:ProofSimpleAction}
\infogen=
X_{i}X_{j}\partial_{X_{i}}\partial_{X_{j}}
+
X_{i}Y_{j}\partial_{Y_{i}}\partial_{X_{j}}
+
Y_{i}X_{j}\partial_{X_{i}}\partial_{Y_{j}}
+
Y_{i}Y_{j}\partial_{Y_{i}}\partial_{Y_{j}}
.
\end{gather}
\end{Definition}

\begin{Lemma}\label{L:ProofSimpleEndAlg}
The assignment $\infbgen\mapsto\infogen$
endows $\verma{\oplus\bsym}$ with the structure
of a $\infbraid$ module.
This $\infbraid$ action stabilizes 
$(\verma{\oplus\bsym})_{\bsym[d]}$ for all 
$\bsym[d]\in\Z^{n}$.
\end{Lemma}

\begin{proof}
A direct calculation, see also \cite[Theorem 2.1]{LeZh-strongly-multiplicity-free}.
\end{proof}

\begin{Remark}\label{R:ProofSimpleOperators}
The $\infbraid$ action on $\verma{\oplus\bsym}$
using \autoref{Eq:ProofSimpleAction} factors through 
an $\infbraid$ action on $\ugl[2]^{\otimes n}$, see 
\cite[Section 2]{Ko-cft-topology}. This works 
as follows. Let $B=\{E,F,L_{1},L_{2}\}$ be the 
usual basis of $\gl[2]$, and $E^{\star}=F$, $F^{\star}=E$, $L_{i}^{\star}=L_{i}$.
Then define 
\emph{Casimir-type elements} 
by $\casimir[{ij}]=\sum_{b\in B}
1^{\otimes i-1} 
\otimes b\otimes
1^{\otimes j-i-1}
\otimes b^{\star}\otimes 1^{\otimes n-j}\in\ugl[2]^{\otimes n}$ 
for $1\leq i<j\leq k$. Then $\infbgen\mapsto\casimir[{ij}]$ 
defines an $\infbraid$ action that factors 
the $\infbraid$ action given by $\infbgen\mapsto\infogen$.
\end{Remark}

We have $\verma{\oplus\bsym}\cong\medoplus_{\bsym[d]\in\Z^{n}}
(\verma{\oplus\bsym})_{\bsym[d]}$ 
as $\infbraid$ modules, by construction. But we in the end only 
need a fixed arbitrary direct summand $(\verma{\oplus\bsym})_{\bsym[d]}\cong\verma{\para_{1}+d_{1}}
\otimes\dots\otimes\verma{\para_{n}+d_{n}}$.
Let $\ealg$ denote the image of the $\infbraid$ action from 
\autoref{L:ProofSimpleEndAlg} restricted to $(\verma{\oplus\bsym})_{\bsym[d]}$.
Note that 
$\ealg\subset\End_{\K}\big((\verma{\oplus\bsym})_{\bsym[d]}\big)$ but we will need 
the following stronger statement.

\begin{Lemma}\label{L:ProofSimpleEndAlgUsingHowe}
Assume that the parameters are admissible.
We have $\ealg\subset\End_{\ugl[2]}\big((\verma{\oplus\bsym})_{\bsym[d]}\big)$ and 
$\ealg\dgen\End_{\ugl[2]}\big((\verma{\oplus\bsym})_{\bsym[d]}\big)$.
\end{Lemma}

\begin{proof}
The proof splits into several parts.

\textit{Containment.}
$\ealg\subset\End_{\ugl[2]}\big((\verma{\oplus\bsym})_{\bsym[d]}\big)$ 
follows from \autoref{Eq:ProofSimpleAction} via a direct computation.

\textit{Applying Verma Howe duality.}
From \autoref{T:DualityHowe} we get
commuting actions of $\ugl[2]$ and $\ugl[n]$ on 
$\verma{\oplus\bsym}$ and the $\ugl[2]$-$\ugl[n]$ bimodule 
decomposition
\begin{gather*}
(\verma{\oplus\bsym})_{\bsym[d]}
\cong
\medoplus_{t\in\N}
\verma{\sums{\para_{n}}+\countingpara-t,t}
\otimes(\dense{\countingpara-t,t})_{\bsym[d]}
,
\end{gather*}
with $\verma{\sums{\para_{n}}+\countingpara-t,t}$ and $\dense{\countingpara-t,t}$ being simple.
It follows that
\begin{gather*}
\End_{\ugl[2]}\big((\verma{\oplus\bsym})_{\bsym[d]}\big)
\cong
\End_{\K}\big(
\medoplus_{t\in\N}(\dense{\countingpara-t,t})_{\bsym[d]}\big).
\end{gather*}
It 
remains to show that all endomorphisms of 
$(\dense{\countingpara-t,t})_{\bsym[d]}$ come from $\ealg$ 
in the dense sense.

\textit{The Casimir subalgebra.}
To this end, let $\casalg\subset
\End_{\K}\big((\verma{\oplus\bsym}_{\bsym[d]})\big)$, by definition, 
be the $\K$ algebra generated by
$\langle e_{ij}e_{ji}+e_{ji}e_{ij},e_{kk}|1\leq i<j\leq n,1\leq k\leq n\rangle$ with the endomorphisms as in \autoref{Eq:ProofDiffOpTwo}.
It is important to observe that the images of the Casimir operators 
$\casimir[k]$ from \autoref{Eq:ProofCasimirElement} are in 
$\casalg$. Moreover, 
note that $\casalg\subset\End_{\ugl[2]}\big((\verma{\oplus\bsym})_{\bsym[d]}\big)$ by 
\autoref{T:DualityHowe} and the elements of $\casalg$ 
are homogeneous, and hence, 
$\casalg$ acts on $(\dense{\countingpara-t,t})_{\bsym[d]}$.

\textit{Simplicity.}
We aim to show that $(\dense{\countingpara-t,t})_{\bsym[d]}$ 
is simple as a $\casalg$ module. For this we use an analog
of \cite[Theorem 6.1 and Remark 6.2]{MiToLa-casimir-braids}, 
where the main observation is that the Casimir operators 
$\casimir[k]$ from \autoref{Eq:ProofCasimirElement} 
have a joint simple spectrum on $(\dense{\countingpara-t,t})_{\bsym[d]}$ 
with diagonal basis given by the GT vectors. 
This follows from the proof of \autoref{L:ProofCasimir}.

Using the GT formulas from \autoref{D:DualityGTModule} it is 
not hard to see that the action graph of the $\casimir[k]$ 
on $(\dense{\countingpara-t,t})_{\bsym[d]}$ with vertices given 
by the relevant GT vectors is strongly connected. Hence, 
as soon as a nonzero $\casalg$ submodule $\module[M]\subset(\dense{\countingpara-t,t})_{\bsym[d]}$ 
contains at least one GT vector we have 
$\module[M]=(\dense{\countingpara-t,t})_{\bsym[d]}$.
Finally, since the spectrum of the $\casimir[k]$ 
is simple with diagonal basis given by the GT vectors by the above, 
every such $\module[M]\subset(\dense{\countingpara-t,t})_{\bsym[d]}$ 
contains indeed a GT vector.

\textit{Wrap-up.}
It follows that $\casalg$ generates the whole of $\End_{\K}\big((\dense{\countingpara-t,t})_{\bsym[d]}\big)$. 
A calculation then verifies that 
$e_{ij}e_{ji}=\infogen+e_{ii}$ and $e_{ii}$ acts as a scalar,
which completes the proof since 
the elements $\infogen$ generate $\ealg$, by definition, and we then have
\begin{gather*}
\ealg\dgen
\End_{\K}\big(
\medoplus_{t\in\N}(\dense{\countingpara-t,t})_{\bsym[d]}\big)
\cong
\End_{\ugl[2]}\big((\verma{\oplus\bsym})_{\bsym[d]}\big)
\end{gather*}
which gives denseness.
\end{proof}

%%%%%%%%%%%%%%%%%%%%%%%%%%%%%%%%%%%%%%%%%

\subsection{The quantum case}\label{SS:ProofSimpleQuantum}

%%%%%%%%%%%%%%%%%%%%%%%%%%%%%%%%%%%%%%%%%

Instead of infinitesimal braids we come back to the pure braid groups.
To this end, recall that 
$(\qverma{\oplus\bsym})_{\bsym[d]}$ is a $\ibraid{\bsym[\rho]}$ module for 
all $\bsym[\rho]\leq\bsym$
by \autoref{L:LKBAction}. 
In particular, $(\qverma{\oplus\bsym})_{\bsym[d]}$ is a $\pbraid$ module.
Let $\eqalg$ be the image of this action, 
and similarly for admissible parameter we let $\evalg$ be the respective image.

\begin{Lemma}\label{L:ProofSimpleDefOverAForm}
The $\Kvpara$ algebra $\evalg$ contains an $\apara$ subalgebra 
$\evalg$ whose classical specialization contains $\ealg$
and whose quantum specialization is $\eqalg$.
\end{Lemma}

\begin{proof}
First note that all appearing scalars are in $\apara\subset\Kvpara$, 
so the construction in \autoref{SS:LKBSimple} works {\ver} over $\apara$.

The statement about the quantum specialization is then clear since 
this is how the quantum specialization is defined.

For the classical specialization the same argument as 
\cite[Proof of Theorem 7.5]{LeZh-strongly-multiplicity-free} works.
\end{proof}

\begin{Lemma}\label{L:ProofSimplePureSurjective}
Assume that the parameters are admissible.\!\!
Then we have $\eqalg\subset\End_{\uqgl[2]}\big((\qverma{\oplus\bsym})_{\bsym[d]}\big)$ and 
$\eqalg\dgen\End_{\uqgl[2]}\big((\qverma{\oplus\bsym})_{\bsym[d]}\big)$.
\end{Lemma}

\begin{proof}
$\eqalg\subset\End_{\uqgl[2]}
\big((\qverma{\oplus\bsym})_{\bsym[d]}\big)$
follows from \autoref{L:LKBAction} 
and	the definition of the $R$ matrices.

For the second statement recall that $\averma{\oplus\bsymv}$ 
is flat and specialize to 
$\verma{\oplus\bsym}$ classically and to $\qverma{\oplus\bsym}$ 
in the quantum case, see \autoref{L:ProofFlat}. Thus, on the side of the 
endomorphism algebra we can change between the classical and 
the quantum case. Moreover, \autoref{L:ProofSimpleDefOverAForm} 
shows that $\eqalg$ is at least as big as $\ealg$.
Taking both together and using \autoref{L:ProofSimpleEndAlgUsingHowe} 
implies then the claim.
\end{proof}

We are ready for the final proof of this paper:

\begin{proof}[Proof of \autoref{T:LKBSimple}]
It is enough to consider $\ibraid{\bsym[\rho]}=\pbraid$, 
so we restrict to this case.	

Note that $(\qverma{\oplus\bsym})_{\bsym[d]}$ is a 
$\uqgl[2]$-$\pbraid$ bimodule by \autoref{L:LKBAction}, 
and moreover \autoref{L:ProofSimplePureSurjective} shows that 
$\eqalg$ densely-generates the centralizer of $\uqgl[2]$ 
on $(\qverma{\oplus\bsym})_{\bsym[d]}$. Having this and the 
usual statements about simple modules of centralizers as {\eg} in
\cite[Theorem 4.2.1]{GoWa-representations-invariants}, 
it remains to argue that the LKB representations are 
$\pbraid$ submodules of some
$(\dense{\countingpara-t,t})_{\bsym[d]}$.
(Note hereby that the LKB story is finite dimensional 
and densely generates turns into generates, and hence, 
\cite[Theorem 4.2.1]{GoWa-representations-invariants} applies.)

To see this we note that, as in the proof of the classical version of \autoref{T:DualityHowe},
$\dense{\countingpara-t,t}$ consists of highest weight vectors 
for the $\uqgl[2]$ action, so the condition on $\qlkb{{n,l}}$ 
to be annihilated by $e$ holds. Moreover, by \autoref{Eq:ProofDiffOpQuantum} we get that $K$ acts on 
$\dense{\countingpara-t,t}$ as the scalar $\qpar^{\sums{\para_{n}}+\countingpara}$. 
In particular, for $\countingpara=-2l$ we get 
$\ker(k-\qpar^{\sums{\para_{n}}-2l})\subset
\dense{\countingpara-t,t}$.
Hence, the LKB representations are 
$\pbraid$ submodules of some
$(\dense{\countingpara-t,t})_{\bsym[d]}$ 
(for some $\bsym[d]\in\Z^{n}$ depending on $l$) as desired.
\end{proof}


\begin{thebibliography}{NUW96}
	
	\bibitem[APW91]{AnPoWe-representation-qalgebras}
	H.H. Andersen, P.~Polo, and K.X. Wen.
	\newblock Representations of quantum algebras.
	\newblock {\em Invent. Math.}, 104(1):1--59, 1991.
	\newblock \href {https://doi.org/10.1007/BF01245066}
	{\path{doi:10.1007/BF01245066}}.
	
	\bibitem[AST18]{AnStTu-cellular-tilting}
	H.H. Andersen, C.~Stroppel, and D.~Tubbenhauer.
	\newblock Cellular structures using {$\mathrm{U}_q$}-tilting modules.
	\newblock {\em Pacific J. Math.}, 292(1):21--59, 2018.
	\newblock URL: \url{https://arxiv.org/abs/1503.00224}, \href
	{https://doi.org/10.2140/pjm.2018.292.21}
	{\path{doi:10.2140/pjm.2018.292.21}}.
	
	\bibitem[AST17]{AnStTu-semisimple-tilting}
	H.H. Andersen, C.~Stroppel, and D.~Tubbenhauer.
	\newblock Semisimplicity of {H}ecke and (walled) {B}rauer algebras.
	\newblock {\em J. Aust. Math. Soc.}, 103(1):1--44, 2017.
	\newblock URL: \url{https://arxiv.org/abs/1507.07676}, \href
	{https://doi.org/10.1017/S1446788716000392}
	{\path{doi:10.1017/S1446788716000392}}.
	
	\bibitem[Big01]{Bi-linear-artin}
	S.J. Bigelow.
	\newblock Braid groups are linear.
	\newblock {\em J. Amer. Math. Soc.}, 14(2):471--486, 2001.
	\newblock URL: \url{https://arxiv.org/abs/math/0005038}, \href
	{https://doi.org/10.1090/S0894-0347-00-00361-1}
	{\path{doi:10.1090/S0894-0347-00-00361-1}}.
	
	\bibitem[Bra37]{Br-brauer-algebra-original}
	R.~Brauer.
	\newblock On algebras which are connected with the semisimple continuous
	groups.
	\newblock {\em Ann. of Math. (2)}, 38(4):857--872, 1937.
	\newblock \href {https://doi.org/10.2307/1968843} {\path{doi:10.2307/1968843}}.
	
	\bibitem[BDK20]{BrDaKu-quantum-type-q-webs}
	G.C. Brown, N.J. Davidson, and J.R. Kujawa.
	\newblock Quantum webs of type Q.
	\newblock 2020.
	\newblock URL: \url{https://arxiv.org/abs/2001.00663}.
	
	\bibitem[CK18]{CaKa-q-satake-sln}
	S.~Cautis and J.~Kamnitzer.
	\newblock Quantum {K}-theoretic geometric {S}atake: the {$\mathrm{SL}_n$} case.
	\newblock {\em Compos. Math.}, 154(2):275--327, 2018.
	\newblock URL: \url{https://arxiv.org/abs/1509.00112}, \href
	{https://doi.org/10.1112/S0010437X17007564}
	{\path{doi:10.1112/S0010437X17007564}}.
	
	\bibitem[CKM14]{CaKaMo-webs-skew-howe}
	S.~Cautis, J.~Kamnitzer, and S.~Morrison.
	\newblock Webs and quantum skew {H}owe duality.
	\newblock {\em Math. Ann.}, 360(1-2):351--390, 2014.
	\newblock URL: \url{https://arxiv.org/abs/1210.6437}, \href
	{https://doi.org/10.1007/s00208-013-0984-4}
	{\path{doi:10.1007/s00208-013-0984-4}}.
	
	\bibitem[CW20]{ChWa-quantum-q-howe}
	Z.~Chang and Y.~Wang.
	\newblock Howe duality for quantum queer superalgebras.
	\newblock {\em J. Algebra}, 547:358--378, 2020.
	\newblock URL: \url{https://arxiv.org/abs/1805.04809}, \href
	{https://doi.org/10.1016/j.jalgebra.2019.11.023}
	{\path{doi:10.1016/j.jalgebra.2019.11.023}}.
	
	\bibitem[CW12]{ChWa-dualities-super}
	S.-J. Cheng and W.~Wang.
	\newblock {\em Dualities and representations of {L}ie superalgebras}, volume
	144 of {\em Graduate Studies in Mathematics}.
	\newblock American Mathematical Society, Providence, RI, 2012.
	\newblock \href {https://doi.org/10.1090/gsm/144} {\path{doi:10.1090/gsm/144}}.
	
	\bibitem[DPS98]{DuPaSc-schur-weyl-tilting}
	J.~Du, B.~Parshall, and L.~Scott.
	\newblock Quantum {W}eyl reciprocity and tilting modules.
	\newblock {\em Comm. Math. Phys.}, 195(2):321--352, 1998.
	\newblock \href {https://doi.org/10.1007/s002200050392}
	{\path{doi:10.1007/s002200050392}}.
	
	\bibitem[DN21]{DuNa-categorical-lkb}
	B. Dupont, G. Naisse.
	\newblock Categorification of infinite-dimensional $\mathfrak{sl}_{2}$-modules and braid group 2-actions {I}: tensor products.
	\newblock 2021.
	\newblock URL: \url{https://arxiv.org/abs/2103.14760}.
	
	\bibitem[ES18]{EhSt-nw-algebras-howe}
	M.~Ehrig and C.~Stroppel.
	\newblock Nazarov--{W}enzl algebras, coideal subalgebras and categorified skew
	{H}owe duality.
	\newblock {\em Adv. Math.}, 331:58--142, 2018.
	\newblock URL: \url{https://arxiv.org/abs/1310.1972}, \href
	{https://doi.org/10.1016/j.aim.2018.01.013}
	{\path{doi:10.1016/j.aim.2018.01.013}}.
	
	\bibitem[Fer90]{Fe-weight-modules}
	S.~L.~Fernando.
	\newblock Lie algebra modules with finite-dimensional weight spaces. {I}.
	\newblock {\em Trans. Amer. Math. Soc.}, 322(2):757--781, 1990.
	\newblock \href {https://doi.org/10.2307/2001724} {\path{doi:10.2307/2001724}}.
	
	\bibitem[For96]{Fo-braid-group-reps}
	E.~Formanek.
	\newblock Braid group representations of low degree.
	\newblock {\em Proc. London Math. Soc.}, (3) 73 (1996), no. 2, 279--322.
	\newblock \href {https://doi.org/10.1112/plms/s3-73.2.279} {\path{doi:10.1112/plms/s3-73.2.279}}.
	
	\bibitem[Fut87]{Fu-weight-modules}
	V.M.~Futorny.
	\newblock The weight representations of semisimple finite dimensional Lie algebras.
	\newblock 1987.
	\newblock Ph.D. thesis, Kiev University.
	
	\bibitem[GW09]{GoWa-representations-invariants}
	R.~Goodman and N.R.~Wallach.
	\newblock {\em Symmetry, representations, and invariants}, volume 255 of {\em
		Graduate Texts in Mathematics}.
	\newblock Springer, Dordrecht, 2009.
	\newblock \href {https://doi.org/10.1007/978-0-387-79852-3}
	{\path{doi:10.1007/978-0-387-79852-3}}.
	
	\bibitem[HOL02]{HaOlLa-handlebodies}
	R.~H{\"a}ring-Oldenburg and S.~Lambropoulou.
	\newblock Knot theory in handlebodies.
	\newblock {\em J. Knot Theory Ramifications}, 11(6):921--943, 2002.
	\newblock Knots 2000 Korea, Vol. 3 (Yongpyong).
	\newblock URL: \url{https://arxiv.org/abs/math/0405502}, \href
	{https://doi.org/10.1142/S0218216502002050}
	{\path{doi:10.1142/S0218216502002050}}.
	
	\bibitem[How95]{Ho-perspectives-invariant-theory}
	R.~Howe.
	\newblock Perspectives on invariant theory: {S}chur duality, multiplicity-free
	actions and beyond.
	\newblock In {\em The {S}chur lectures (1992) ({T}el {A}viv)}, volume~8 of {\em
		Israel Math. Conf. Proc.}, pages 1--182. Bar-Ilan Univ., Ramat Gan, 1995.
	
	\bibitem[How89]{Ho-remarks-invariant-theory}
	R.~Howe.
	\newblock Remarks on classical invariant theory.
	\newblock {\em Trans. Amer. Math. Soc.}, 313(2):539--570, 1989.
	\newblock \href {https://doi.org/10.2307/2001418} {\path{doi:10.2307/2001418}}.
	
	\bibitem[Hum08]{Hu-cat-o}
	J.E.~Humphreys.
	\newblock {\em Representations of semisimple {L}ie algebras in the {BGG} category {$\mathcal{O}$}}, {\em Graduate Studies in Mathematics}, 94. 
	\newblock American Mathematical Society, Providence, RI, 2008. xvi+289 pp.
	\newblock \href {https://doi.org/10.1090/gsm/094} {\path{doi:10.1090/gsm/094}}.
	
	\bibitem[ILZ21]{IoLeZh-verma-schur-weyl}
	K.~Iohara, G.~Lehrer, and R.~Zhang.
	\newblock Equivalence of a tangle category and a category of infinite
	dimensional {$\mathrm{U}_q(\mathfrak{sl}_2)$}-modules.
	\newblock {\em Represent. Theory}, 25:265--299, 2021.
	\newblock URL: \url{https://arxiv.org/abs/1811.01325}, \href
	{https://doi.org/10.1090/ert/568} {\path{doi:10.1090/ert/568}}.
	
	\bibitem[JK11]{JaKe-verma-lkb}
	C.~Jackson and T.~Kerler.
	\newblock The {L}awrence--{K}rammer--{B}igelow representations of the braid
	groups via {$U_q(\mathfrak{sl}_2)$}.
	\newblock {\em Adv. Math.}, 228(3):1689--1717, 2011.
	\newblock URL: \url{https://arxiv.org/abs/0912.2114}, \href
	{https://doi.org/10.1016/j.aim.2011.06.027}
	{\path{doi:10.1016/j.aim.2011.06.027}}.
	
	\bibitem[Jan96]{Ja-lectures-qgroups}
	J.C. Jantzen.
	\newblock {\em Lectures on quantum groups}, volume~6 of {\em Graduate Studies
		in Mathematics}.
	\newblock American Mathematical Society, Providence, RI, 1996.
	
	\bibitem[Jim86]{Ji-q-schur-weyl}
	M.~Jimbo.
	\newblock A {$q$}-analogue of {$U(\mathfrak{gl}(N+1))$}, {H}ecke algebra, and
	the {Y}ang--{B}axter equation.
	\newblock {\em Lett. Math. Phys.}, 11(3):247--252, 1986.
	\newblock \href {https://doi.org/10.1007/BF00400222}
	{\path{doi:10.1007/BF00400222}}.
	
	\bibitem[Jon85]{Jo-jones-polynomial}
	V.F.R. Jones.
	\newblock A polynomial invariant for knots via von {N}eumann algebras.
	\newblock {\em Bull. Amer. Math. Soc. (N.S.)}, 12(1):103--111, 1985.
	\newblock \href {https://doi.org/10.1090/S0273-0979-1985-15304-2}
	{\path{doi:10.1090/S0273-0979-1985-15304-2}}.
	
	\bibitem[K{\aa}h10]{Ka-tensor-infinite-dimensional-modules}
	J.~K{\aa}hrstr{\"o}m.
	\newblock Tensoring with infinite-dimensional modules in {$\mathcal{O}_0$}.
	\newblock {\em Algebr. Represent. Theory}, 13(5):561--587, 2010.
	\newblock URL: \url{https://arxiv.org/abs/0708.2218}, \href
	{https://doi.org/10.1007/s10468-009-9137-6}
	{\path{doi:10.1007/s10468-009-9137-6}}.
	
	\bibitem[Kas95]{Ka-quantum-groups}
	C.~Kassel.
	\newblock {\em Quantum groups}, volume 155 of {\em Graduate Texts in
		Mathematics}.
	\newblock Springer-Verlag, New York, 1995.
	\newblock \href {https://doi.org/10.1007/978-1-4612-0783-2}
	{\path{doi:10.1007/978-1-4612-0783-2}}.
	
	\bibitem[Koh02]{Ko-cft-topology}
	T.~Kohno.
	\newblock {\em Conformal field theory and topology}, volume 210 of {\em
		Translations of Mathematical Monographs}.
	\newblock American Mathematical Society, Providence, RI, 2002.
	\newblock Translated from the 1998 Japanese original by the author, Iwanami
	Series in Modern Mathematics.
	\newblock \href {https://doi.org/10.1090/mmono/210}
	{\path{doi:10.1090/mmono/210}}.
	
	\bibitem[Koh12]{Ko-quantum-braid-lkb}
	T.~Kohno.
	\newblock Quantum and homological representations of braid groups.
	\newblock In {\em Configuration spaces}, volume~14 of {\em CRM Series}, pages
	355--372. Ed. Norm., Pisa, 2012.
	\newblock \href {https://doi.org/10.1007/978-88-7642-431-1_16}
	{\path{doi:10.1007/978-88-7642-431-1_16}}.
	
	\bibitem[Kra02]{Kr-linear-artin}
	D.~Krammer.
	\newblock Braid groups are linear.
	\newblock {\em Ann. of Math. (2)}, 155(1):131--156, 2002.
	\newblock URL: \url{https://arxiv.org/abs/math/0405198}, \href
	{https://doi.org/10.2307/3062152} {\path{doi:10.2307/3062152}}.
	
	\bibitem[LNV21]{LaNaVa-tensor-product-blob}
	A.~Lacabanne and G.~Naisse and P.~Vaz.
	\newblock Tensor product categorifications, {V}erma modules and the blob {2}-category.
	\newblock {\em Quantum Topol.}, 12 (2021), no. 4, 705--812.
	\newblock URL: \url{https://arxiv.org/abs/2005.06257}, \href
	{https://doi.org/10.4171/QT/156}
	{\path{doi:10.4171/QT/156}}.
	
	\bibitem[LTV23]{LaTuVa-annular-webs-levi}
	A.~Lacabanne, D.~Tubbenhauer, and P.~Vaz.
	\newblock Annular webs and {L}evi subalgebras.
	\newblock {\em J. Comb. Algebra} 7 (2023), no. 3-4, 283--326.
	\newblock URL: \url{https://arxiv.org/abs/2204.00947}, \href
	{https://doi.org/10.4171/jca/76}
	{\path{doi:10.4171/jca/76}}.
	
	\bibitem[LV21]{LaVa-schur-weyl-ariki-koike}
	A.~Lacabanne and P.~Vaz.
	\newblock {S}chur--{W}eyl duality, {V}erma modules, and row quotients of
	{A}riki--{K}oike algebras.
	\newblock {\em Pacific J. Math.}, 311(1):113--133, 2021.
	\newblock URL: \url{https://arxiv.org/abs/2004.01065}, \href
	{https://doi.org/10.2140/pjm.2021.311.113}
	{\path{doi:10.2140/pjm.2021.311.113}}.
	
	\bibitem[Law90]{La-lkb-rep}
	R.J. Lawrence.
	\newblock Homological representations of the {H}ecke algebra.
	\newblock {\em Comm. Math. Phys.}, 135(1):141--191, 1990.
	\newblock \href {https://doi.org/10.1007/BF02097660}
	{\path{doi:10.1007/BF02097660}}.
	
	\bibitem[LZZ11]{LeZhZh-q-first-fundamental-theorem}
	G.I. Lehrer, H.~Zhang, and R.B. Zhang.
	\newblock A quantum analogue of the first fundamental theorem of classical
	invariant theory.
	\newblock {\em Comm. Math. Phys.}, 301(1):131--174, 2011.
	\newblock URL: \url{https://arxiv.org/abs/0908.1425}, \href
	{https://doi.org/10.1007/s00220-010-1143-3}
	{\path{doi:10.1007/s00220-010-1143-3}}.
	
	\bibitem[LZ06]{LeZh-strongly-multiplicity-free}
	G.I. Lehrer and R.B. Zhang.
	\newblock Strongly multiplicity free modules for {L}ie algebras and quantum
	groups.
	\newblock {\em J. Algebra}, 306(1):138--174, 2006.
	\newblock \href {https://doi.org/10.1016/j.jalgebra.2006.03.043}
	{\path{doi:10.1016/j.jalgebra.2006.03.043}}.
	
	\bibitem[Lus90]{Lu-qgroups-root-of-1}
	G.~Lusztig.
	\newblock Quantum groups at roots of {$1$}.
	\newblock {\em Geom. Dedicata}, 35(1-3):89--113, 1990.
	\newblock \href {https://doi.org/10.1007/BF00147341}
	{\path{doi:10.1007/BF00147341}}.
	
	\bibitem[Mar22]{Ma-homological-colored-lkb}
	J.~Martel.
	\newblock A homological model for {$U_{q}\mathfrak{sl}_{2}$} {V}erma modules and their
	braid representations.
	\newblock  {\em Geom. Topol.} 26 (2022), no. 3, 1225--1289.
	\newblock URL: \url{https://arxiv.org/abs/2002.08785}, \href
	{https://doi.org/10.2140/gt.2022.26.1225}
	{\path{doi:10.2140/gt.2022.26.1225}}.
	
	\bibitem[Mar20]{Ma-colored-lkb}
	J.~Martel.
	\newblock Colored version for {L}awrence representations.
	\newblock 2020.
	\newblock To appear in {\em New York J. Math.}
	\newblock URL: \url{https://arxiv.org/abs/2004.00977}.
	
	\bibitem[Mat00]{Ma-classification-weight-modules}
	O.~Mathieu.
	\newblock Classification of irreducible weight modules.
	\newblock {\em Ann. Inst. Fourier (Grenoble)}, 50(2):537--592, 2000.
	\newblock \href {https://doi.org/10.5802/aif.1765}
	{\path{doi:10.5802/aif.1765}}.
	
	\bibitem[Maz03]{Ma-gt-dense}
	V.~Mazorchuk.
	\newblock Quantum deformation and tableaux realization of simple dense
	{${\mathfrak{gl}}(n,\mathbbm{C})$}-modules.
	\newblock {\em J. Algebra Appl.}, 2(1):1--20, 2003.
	\newblock \href {https://doi.org/10.1142/S0219498803000325}
	{\path{doi:10.1142/S0219498803000325}}.
	
	\bibitem[MTL05]{MiToLa-casimir-braids}
	J.J. Millson and V.~Toledano~Laredo.
	\newblock Casimir operators and monodromy representations of generalised braid
	groups.
	\newblock {\em Transform. Groups}, 10(2):217--254, 2005.
	\newblock URL: \url{https://arxiv.org/abs/math/0305062}, \href
	{https://doi.org/10.1007/s00031-005-1008-6}
	{\path{doi:10.1007/s00031-005-1008-6}}.
	
	\bibitem[NS95]{NoSu-q-symmetric-spaces}
	M. Noumi and T. Sugitani.
	\newblock Quantum symmetric spaces and related {$q$}-orthogonal polynomials.
	\newblock In {\em Group theoretical methods in physics ({T}oyonaka, 1994)},
	pages 28--40. World Sci. Publ., River Edge, NJ, 1995.
	
	\bibitem[NUW96]{NoUmWa-sl2-son-duality}
	M.~Noumi, T.~Umeda, and M.~Wakayama.
	\newblock Dual pairs, spherical harmonics and a {C}apelli identity in quantum
	group theory.
	\newblock {\em Compositio Math.}, 104(3):227--277, 1996.
	
	\bibitem[QS19]{QuSa-mixed-skew-howe}
	H.~Queffelec and A.~Sartori.
	\newblock Mixed quantum skew {H}owe duality and link invariants of type {$A$}.
	\newblock {\em J. Pure Appl. Algebra}, 223(7):2733--2779, 2019.
	\newblock URL: \url{https://arxiv.org/abs/1504.01225}, \href
	{https://doi.org/10.1016/j.jpaa.2018.09.014}
	{\path{doi:10.1016/j.jpaa.2018.09.014}}.
	
	\bibitem[QW24]{QuWe-extremal-projectors-2}
	H.~Queffelec and P.~Wedrich.
	\newblock Extremal weight projectors {II}, $\mathfrak{gl}_{N}$ case.
	\newblock {\em Algebr. Comb.} 7 (2024), no. 1, 187–223.
	\newblock URL: \url{https://arxiv.org/abs/1803.09883}, \href
	{https://doi.org/10.5802/alco.330}
	{\path{doi:10.5802/alco.330}}.
	
	\bibitem[RT21]{RoTu-homflypt-handlebody}
	D.E.V. Rose and D.~Tubbenhauer.
	\newblock {HOMFLYPT} homology for links in handlebodies via type {$\mathsf{A}$}
	{S}oergel bimodules.
	\newblock {\em Quantum Topol.}, 12(2):373--410, 2021.
	\newblock URL: \url{https://arxiv.org/abs/1908.06878}, \href
	{https://doi.org/10.4171/qt/152} {\path{doi:10.4171/qt/152}}.
	
	\bibitem[RT16]{RoTu-symmetric-howe}
	D.E.V. Rose and D.~Tubbenhauer.
	\newblock Symmetric webs, {J}ones--{W}enzl recursions, and {$q$}-{H}owe
	duality.
	\newblock {\em Int. Math. Res. Not. IMRN}, (17):5249--5290, 2016.
	\newblock URL: \url{https://arxiv.org/abs/1501.00915}, \href
	{https://doi.org/10.1093/imrn/rnv302} {\path{doi:10.1093/imrn/rnv302}}.
	
	\bibitem[ST19]{SaTu-bcd-webs}
	A.~Sartori and D.~Tubbenhauer.
	\newblock Webs and {$q$}-{H}owe dualities in types {BCD}.
	\newblock {\em Trans. Amer. Math. Soc.}, 371(10):7387--7431, 2019.
	\newblock URL: \url{https://arxiv.org/abs/1701.02932}, \href
	{https://doi.org/10.1090/tran/7583} {\path{doi:10.1090/tran/7583}}.
	
	\bibitem[Sch01]{Sc-schur-weyl}
	I.~Schur.
	\newblock {\"U}ber eine {K}lasse von {M}atrizen, die sich einer gegebenen
	{M}atrix zuordnen lassen.
	\newblock 1901.
	\newblock Ph.D. thesis, Friedrich-Wilhelms--Universit{\"a}t Berlin.
	\newblock URL: \url{https://gdz.sub.uni-goettingen.de/id/PPN271034092}.
	
	\bibitem[TV23]{TuVa-handlebody}
	D.~Tubbenhauer and P.~Vaz.
	\newblock Handlebody diagram algebras.
	\newblock Rev. Mat. Iberoam. 39 (2023), no. 3, 845--896.
	\newblock URL: \url{https://arxiv.org/abs/2105.07049}, \href
	{https://doi.org/10.4171/RMI/1356} {\path{doi:10.4171/RMI/1356}}.
	
	\bibitem[TVW17]{TuVaWe-super-howe}
	D.~Tubbenhauer, P.~Vaz, and P.~Wedrich.
	\newblock Super {$q$}-{H}owe duality and web categories.
	\newblock {\em Algebr. Geom. Topol.}, 17(6):3703--3749, 2017.
	\newblock URL: \url{https://arxiv.org/abs/1504.05069}, \href
	{https://doi.org/10.2140/agt.2017.17.3703}
	{\path{doi:10.2140/agt.2017.17.3703}}.
	
	\bibitem[Ver98]{Ve-handlebodies}
	V.V. Vershinin.
	\newblock On braid groups in handlebodies.
	\newblock {\em Sibirsk. Mat. Zh.}, 39(4):755--764, i, 1998.
	\newblock \href {https://doi.org/10.1007/BF02673050}
	{\path{doi:10.1007/BF02673050}}.
	
	\bibitem[Wey97]{We-classical-groups}
	H.~Weyl.
	\newblock {\em The classical groups}.
	\newblock Princeton Landmarks in Mathematics. Princeton University Press,
	Princeton, NJ, 1997.
	\newblock Their invariants and representations, Fifteenth printing, Princeton
	Paperbacks.
	
	\bibitem[Zin01]{Zi-lkb-bmw}
	M.G. Zinno.
	\newblock On {K}rammer's representation of the braid group.
	\newblock {\em Math. Ann.}, 321(1):197--211, 2001.
	\newblock URL: \url{https://arxiv.org/abs/math/0002136}, \href
	{https://doi.org/10.1007/PL00004501} {\path{doi:10.1007/PL00004501}}.
	
\end{thebibliography}
\end{document}